\documentclass{elsart}

\usepackage{natbib}
\usepackage{overpic,graphicx}
\usepackage{amssymb,amsfonts,amsmath}
\usepackage{color}

\usepackage{par06}

 \def\Xw{{\mathbf X}}
 \def\Yw{{\mathbf Y}}
 \def\Vw{{\mathbf V}}
 \def\Db{{\bar{D}}}
 \def\Sb{{\bar{S}}}
\def\kb{{\bar{k}}}

\def\Qed#1{\leavevmode\unskip\penalty9999 \hbox{}\nobreak\hfill\quad\hbox{#1}}
\def\defumg#1#2#3#4#5{\newenvironment{#1}{\begin{trivlist}%
        \item[\hskip\labelsep{#5#2}]#3}{#4\end{trivlist}}}
\defumg{proof}{Proof.}{}{\Qed{$\square$}}{\it}
\newdimen\Columnwidth
\begin{document}

\begin{frontmatter}


% use the thanksref command within \title, \author or \address for footnotes
\title{Darboux Cyclides and Webs from Circles}
\author{Helmut Pottmann, Ling Shi and Mikhail Skopenkov}
\address{King Abdullah University of Science and Technology, Thuwal, Saudi Arabia}

\begin{abstract}
Motivated by potential applications in architecture, we study Darboux cyclides. These algebraic surfaces
of order $\le 4$ are a superset of Dupin cyclides and quadrics, and they carry up to six real families of circles. Revisiting the  classical approach to
these surfaces based on the spherical model of 3D M\"obius geometry, we provide computational tools for
the identificiation of circle families on a given cyclide and for the direct design of those. In particular,
we show that certain triples of circle families may be arranged as so-called hexagonal webs, and we
provide a complete classification of all possible hexagonal webs of circles on Darboux cyclides.
\end{abstract}
\end{frontmatter}

\section{Introduction}
%%%%%%%%%%%%%%%%%%%%%%

The development of Computer-Aided Geometric Design has largely been driven by applications. Its origins trace back to the need for computationally efficient design and processing of freeform surfaces in the automotive, airplane and ship industry. Many other applications (industrial, medical, scientific data analysis and visualization) shaped the development of the field. However, architecture --- which has been using CAD since decades --- only very recently attracted the attention of the Geometric Modeling community (see \cite{archgeom:2007}).  The main reason for the growing interest of Geometric Modeling in problems from architecture is the trend towards freeform structures which is very clearly seen in contemporary architecture and led by star architects such as Frank Gehry or Zaha Hadid.

While digital models of architectural freeform surfaces are easily created using standard modeling tools, the actual fabrication and construction of architectural freeform structures remains a challenge. In order to make a freeform design realizable, an optimization process known as {\it rationalization} has to be applied.  This means its decomposition into smaller parts, thereby meeting two competing objectives: feasibility, and consistency with the designer's intentions. Most rationalization techniques replace smooth surfaces, possibly with an additional curve network
(of panel boundaries) on them, by other structures like meshes with special properties. The guiding thought in all considerations is the efficient manufacturing of the surface parts (panels) and their respective necessary supporting/connecting elements. Both simple geometry and repetition of elements contribute to this goal of efficiency (see e.g. \cite{liu2006,pottmann-2007-pm,pottmann-2008-strip,schiftner-2009-cp,eigensatz_paneling_sig_10,k-set-tilable:2010,singh_triEquivalence_sig_10}).

\begin{figure}[htb]
\centering{%
\begin{overpic}[width=.68\textwidth]{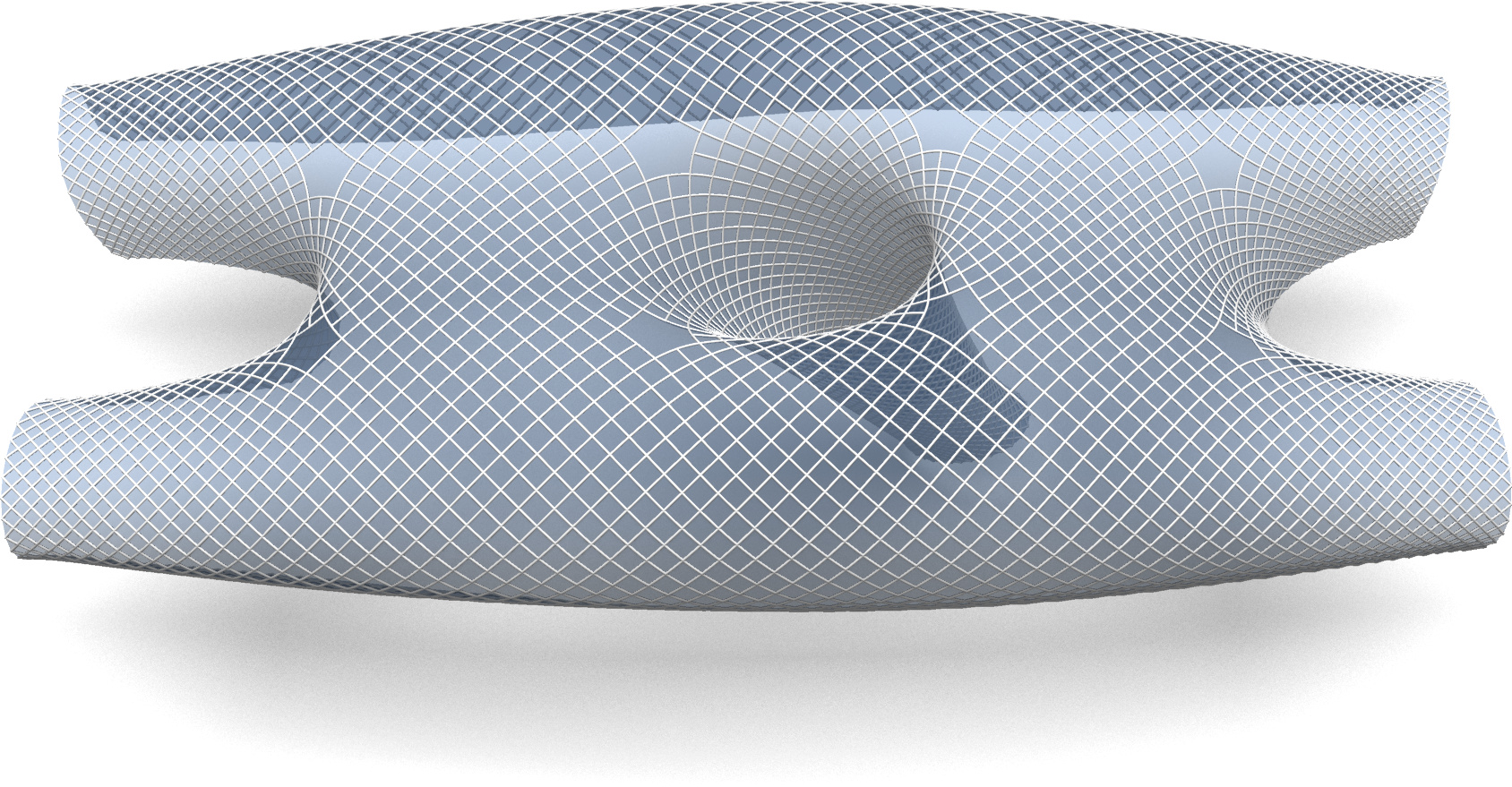}
\end{overpic}
\begin{overpic}[width=.31\textwidth]{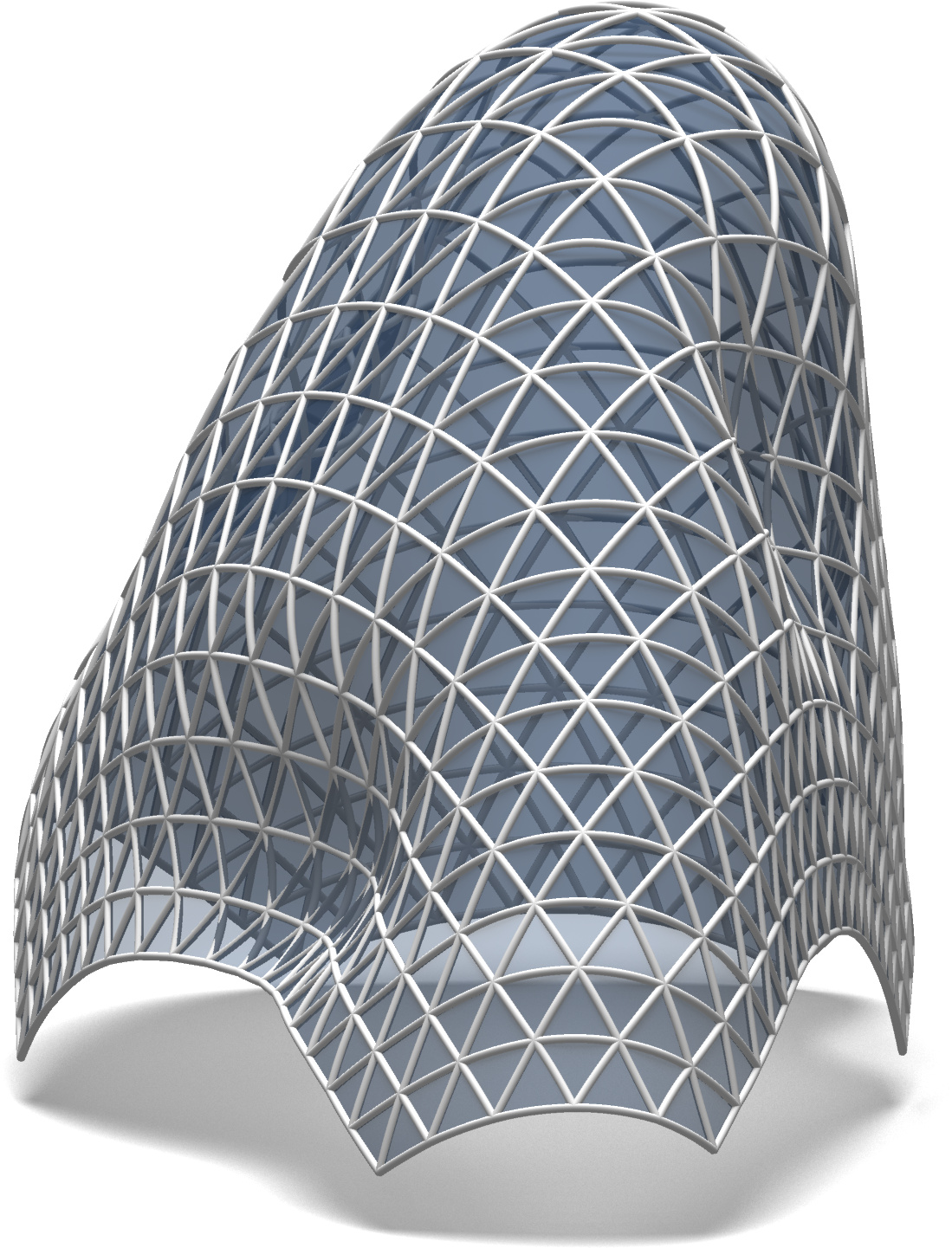}
\end{overpic}
}
 \caption{Circular arc structures according to \cite{bo-2011}. Left: Quadrilateral CAS with an edge angle of 90 degrees covering an architectural design. Right: Triangular CAS (edge angle 60 degrees) on the
Eindhoven Blob by M. Fuksas.}
	\label{fig:cas}
\end{figure}

As a contribution towards rationalization of architectural freeform structures, \cite{bo-2011} have recently
studied so-called {\it circular arc structures}. A circular arc structure (CAS) is a mesh whose edges are realized as circular arcs instead of straight line segments and which possesses congruent nodes with well-defined tangent planes.  Figure~\ref{fig:cas} shows two examples of circular arc structures with quad mesh and triangle mesh combinatorics, respectively. The nodes are congruent and possess edge angles of 90 and 60 degrees, respectively.

A circular arc structure (CAS) of regular quad mesh combinatorics consists of two discrete sets
of arc splines (tangent continuous curves formed by circular arcs) which intersect each other under constant,
in particular right angle. For a triangular CAS, we have three
sets of arc splines meeting in regular 60 degree nodes. Three families of curves which are the
image of a regular triangular grid under a diffeomorphism form a so-called {\it hexagonal web} or {\it 3-web}
 (see \cite{blaschke-bol:1938}).
Thus, a triangular CAS is a special discrete 3-web.

 Now we leave aside the congruence of nodes, but ask for CAS whose generating curves are single
 arcs and not just arc splines. This amounts to the following basic questions:

\begin{enumerate}
\item Which structures can be built from two families of circles (but not necessarily meeting
under a constant angle) and how can we design them?
\item How can we characterize and design 3-webs formed by circles (\emph{C-webs})?
\end{enumerate}

It is conjectured that a surface which carries three
families of circles is a \emph{Darboux cyclide} (defined in Section~\ref{sec:basic});
%or a limit case of it;
a proof is missing so far.
However, it seems to be reasonable to assume that the conjecture is correct; there
is no known surface different from a Darboux cyclide which carries three families of circles.
There are, however, surfaces with two families of circles which are not Darboux cyclides: an example
is furnished by a translational surface with circles as profile curves. To keep the present
paper more focused and since we are mostly interested in webs from circles, we will only study
Darboux cyclides. Knowing that a complete answer to the first question would require additional
efforts, we may return to this topic in another publication.

In the present paper, we will lay the geometric fundamentals for an understanding of Darboux cyclides
(which hopefully will also serve as an introduction to basic concepts of 3D M\"obius geometry)
and then proceed towards design methods for Darboux cyclides, their families of circles and C-webs.
Our main result is a classification of all possible C-webs on Darboux cyclides (Theorem~\ref{thm-main}).

\subsection{Previous work}
%----------------------------------------------------------------------------------------
%
Darboux cyclides are a topic of classical geometry. Main contributions
are due to \cite{kummer:1865} and \cite{darboux:1880}. Probably the most complete discussion
is found in the monograph by  \cite{coolidge:1916}.  After a long period of silence on these
surfaces, geometers again got fascinated by them, especially since they can carry up to six
 families of real circles (see \cite{blum:1980,takeuchi:2000}). Related studies are those which show under certain
 assumptions that a surface with a higher number of circle families must be a sphere (see \cite{takeuchi:1985,takeuchi:1987}).

In {\em geometric modeling}, Dupin cyclides and various extensions attracted a lot of attention (see the survey by \cite{degen:2002}), but to our best knowledge, Darboux cyclides have not yet been studied from the geometric design perspective, with the noteable exception of a very recent paper by \cite{krasauskas:2011}, whose relation to cyclides has been realized in a discussion of our research with Rimas Krasauskas.

The problem of determining {\em all hexagonal webs from circles} in the plane (or equivalently, on the sphere) turned
out to be very difficult. A special class of webs from circles, initially found
by \cite{volk-1929-cweb} by analytical methods, has soon been realized as a simple projection of webs from
straight lines in the plane onto the sphere (\cite{strubecker-1932}).
 An elegant
construction of a nontrivial class of webs from circles by \cite{wunderlich-38} will turn out to be
an important ingredient in our study of C-webs on Darboux cyclides.
Recently,  \cite{shelekhov-2007-cw} could classify all webs formed by pencils of circles.  It is still unknown whether the so far described instances of webs from circles are the only possible ones.

\noindent {\bf Overview.} Our paper is organized as follows.
Section \ref{sec:basic} studies Darboux cyclides by employing the spherical model of M\"obius geometry and
mostly methods of projective geometry (pencils of quadrics) to get insight into the circle families
they contain.  Based on these insights, we are able to
 devise design methods for Darboux cyclides (Section \ref{sec:design}) and to characterize all C-webs on them
 in Section  \ref{sec:cwebs}.
Finally, we address future research directions and open problems in Section  \ref{sec:future}.
%\bigskip

%%%%%%%%%%%%%%%%%%%%%%%%%%%%%%%%%%%%%%%%%%%%%%%%%%%%%%%%%%%%%%%%%%%%%%%%%%%%%%
\section{Darboux cyclides and their families of circles} \label{sec:basic}
%%%%%%%%%%%%%%%%%%%%%%%%%%%%%%%%%%%%%%%%%%%%%%%%%%%%%%%%%%%%%%%%%%%%%%%%%%%%

 We start this section with the definition of a Darboux cyclide. In Section~\ref{ssec:basic} we introduce the spherical model of M\"obius geometry and in Section~\ref{ssec:pencil} we show that a Darboux cyclide is represented by the carrier of a pencil of quadrics in this model. Using this representation, we describe families of circles on a Darboux cyclide in Section~\ref{ssec:circles} and determine the possible number of such families in Section~\ref{ssec:realcircles}. In Section~\ref{ssec:polardual} we apply polarity in the spherical model to set up some background for design of Darboux cyclides in the next section.

%\bigskip
A {\em Darboux cyclide} is a surface whose equation in a Cartesian coordinate system has the form
\be \label{D}
D:\ \lambda(x^2+y^2+z^2)^2 + (x^2+y^2+z^2)L(x,y,z) + Q(x,y,z)= 0,
\ee
with a constant $\lambda$, a polynomial
$L=\mu x+\nu y+\kappa z$, and a %quadratic
polynomial $Q$ of degree at most $2$, not vanishing simultaneously. Hereafter the coefficients and the coordinates are real numbers unless otherwise explicitly indicated.

If the left-hand side of equation~\eqref{D} factors into non-constant polynomials in $x,y,z$ with (possibly complex) coefficients, then the Darboux cyclide is called \emph{reducible}. A reducible cyclide either splits into a union of spheres/planes or degenerates to a curve in $\mathbb{R}^3$.
%If $\lambda=L=0$ and $Q=a(x^2+y^2+z^2)+bx+cy+dz+e$ in equation~\eqref{D}
A cyclide of the form $a(x^2+y^2+z^2)+bx+cy+dz+e=0$
is itself a sphere/plane
%(possibly imaginary or degenerating to a point/plane)
and is referred as \emph{trivial} hereafter. Thus, in the following we focus only on irreducible and nontrivial cyclides.
%A cyclide is \emph{irreducible}, if the left-hand side of equation~(\ref{D}) does not factor into nonconstant polynomials in $x,y,z$ with complex coefficients. Hereafter we focus only on irreducible cyclides.

The
intersection curve of a cyclide $D$ with a sphere $S: \xw^2=L_1(\xw)$, where
 $\xw:=(x,y,z)$ and $L_1$ is a %linear
polynomial of degree at most $1$, lies on the quadric
$$ Q_1:\ \lambda L_1(\xw)^2+L_1(\xw)L(\xw)+Q(\xw)=0,$$
and thus it is in general an algebraic curve $c$ of order at most $4$, which is common
to all quadrics of the pencil spanned by $S$ and $Q_1$.
Hence, if the sphere $S$ touches $D$ at two points, the curve $c$ has two double points and thus it is formed by two circles (or a doubly counted circle $c$,
in which case the sphere touches
$D$ along $c$). This already indicates the
presence of circles on cyclides, which is better studied by a slightly more
advanced method as follows.

\subsection{Some basic facts on the spherical model of M\"obius geometry} \label{ssec:basic}
%--------------------------------------------------------------------------
%
We map points $\xw=(x,y,z)\in \R^3$ to points $\bar{\xw}=(x_1,\ldots,x_4)\in \R^4$ on the sphere $\Sigma: x_1^2+\ldots+x_4^2=1$
via the inverse stereographic projection. The projection center shall be $Z=(0,0,0,1)$ and we
embed $3$-space into $4$-space via $\xw \mapsto (x,y,z,0) \in \R^4$. Corresponding points $\xw$ and $\bar{\xw}$ lie on straight lines through $Z$. The mappings between
$\xw$ and $\bar{\xw}$ in both directions are described by the formulae
 \be \label{stereo1}
 \bar{\xw}=(x_1,\ldots,x_4)={1 \over \xw^2+1} (2x,2y,2z, \xw^2-1),
 \ee
 \be \label{stereo2}
 {\xw}=(x,y,z)={1 \over 1-x_4 } (x_1,x_2,x_3).
 \ee
 %
 %In this way, points of {\em M\"obius geometry} in $3$-space correspond to points of the sphere $\Sigma$; the ideal point $\infty$ of M\"obius geometry corresponds to the projection center $Z$.
By definition, \emph{the set of points in} {\em M\"obius geometry} is Euclidean $3$-space with one \emph{ideal point} $\infty$ added;
 %the inverse stereographic projection takes
 the ideal point $\infty$ corresponds to the projection center $Z$.

\begin{figure}[htb]
\centering{%\vspace{4cm}
\includegraphics
%[width=10.92cm
%width=0.8\textwidth
%]
{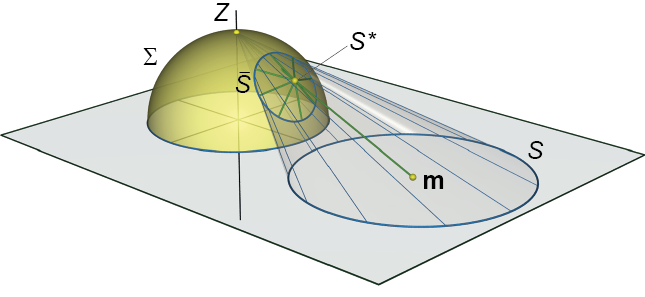}
}
 \caption{The spherical model of planar M\"obius geometry and stereographic projection: any  planar section $\bar S$ of the sphere $\Sigma$ corresponds to an M-circle $S$ and the extended stereographic projection maps the pole $S^*$ of the plane to the center $\mw$ of the circle. The projection center $Z \in \Sigma$ corresponds to the ideal point $\infty$ and circles through $Z$ correspond to straight lines in the plane.}
	\label{fig:stereographic}
\end{figure}

 Unfortunately, it is hard to illustrate geometric facts involving a 4D space. Thus, to support the reader's
 imagination, Fig.~\ref{fig:stereographic} illustrates the lower dimensional counterpart of the stereographic
 projection of a sphere in $3$-space onto a plane. What we are going to say below %have said above
 about %M-
 spheres, is true for %M-
 circles in the lower dimensional counterpart.

 A {\em M\"obius sphere} (or \emph{M-sphere}) $S$ is the set (strictly speaking, in a complex $3$-space) given by the equation $\lambda \xw^2 +ax+by+cz+d=0$, where $\lambda, a, b, c, d$ do not vanish simultaneously. By definition, $\infty\in S$ if and only if $\lambda=0$. A M\"obius sphere is either a plane ($\lambda=0$), or a Euclidean sphere (\emph{real} sphere henceforth), or a point (\emph{null} sphere), or does not have real points at all (\emph{imaginary} sphere).

 One can work with imaginary spheres without appealing to complex numbers as follows. The equation of a nonplanar M-sphere can be rewritten in the form $(\xw-\mw)^2=r^2$, where $\mw\in\mathbb{R}^3$ is the \emph{midpoint} and $r^2\in\mathbb{R}$ is the \emph{square of the radius}. In the case of imaginary sphere $S$ the number $r^2$ is negative; we can write $r^2=-\rho^2$ for some $\rho>0$. Define the {\em real representer} $S^r$ of the sphere $S$ as the sphere with center $\mw$ and radius $\rho$. In the following, all statements involving imaginary spheres can be easily converted to statements on their real representers using Remark~\ref{rem:complex-sphere} below.

 %{\em M\"obius spheres} $S$ (M-spheres), $\lambda \xw^2 +ax+by+cz+d=0$, are either Euclidean spheres ($\lambda \ne 0$) or planes ($\lambda = 0$). Without further restrictions on $(\lambda,a,b,c,d)$, a sphere may also be complex or have radius $0$ (\emph{null sphere}); in our discussions, complex spheres will be of minor importance.

 M-spheres %are mapped via inverse stereographic projection
 correspond to intersections $\bar{S}$ of $\Sigma$ with hyperplanes
 $$ H_S:\ \ ax_1+bx_2+cx_3 + (\lambda - d)x_4 + \lambda+d=0;$$
 of course, the $2$-surfaces $\bar{S}$ are M-spheres, possibly imaginary or null.
 If the M-sphere $S$ is a plane, it contains the ideal point, i.e., $H_S$ and $\bar{S}$ pass through $Z$. Often one uses the polarity with respect to $\Sigma$. Using homogeneous coordinates $X:=(X_1,\ldots,X_5)$ in the projective
 extension $P^4$ of the $4$-space,
 with $x_i=X_i/X_5, \ i=1,\ldots,4$,  the pole $S^*$ of $H_S$ is given by
 $$ S^*=(a,b,c,\lambda-d,-\lambda-d).$$
 The coordinates of $S^*$ are called {\em pentaspherical coordinates} of the M-sphere $S$.
 For an M-sphere $S$ with center $\mw=(m_x,m_y,m_z)$ and squared radius $r^2$, we have
 \be  \lambda=1,\ a=-2m_x, b=-2m_y, c=-2m_z,\  d=\mw^2-r^2. \label{penta} \ee
 The stereographic projection is strictly speaking a mapping from points of the sphere $\Sigma \subset \R^4$
 to $\R^3 \cup \{\infty\}$. However, it is useful to also use the underlying central projection with center $Z$ and image space $x_4=0$ for our considerations, even if we do not apply it only to points of $\Sigma$; we will call this  projection the
 {\em extended stereographic  projection} henceforth. From (\ref{penta}) we can easily verify the well-known fact that the extended stereographic projection  maps the pole $S^*$ onto the midpoint of $S$ (see also Fig.~\ref{fig:stereographic}).

 %% Figure stereographic was here

 An {\em M-circle} $k$ of 3D M\"obius geometry is the intersection of two M-spheres (strictly speaking, in the complex $3$-space) and thus it is either a circle, %(possibly complex or null)
 a straight line, a single point, or the empty set from a Euclidean perspective. Formal definitions of \emph{real}, \emph{null}, and \emph{imaginary} circles are left to the reader. An M-circle is mapped
 %It corresponds
 to the intersection $\bar{k}$ of $\Sigma$ with a 2D plane, i.e., to an M-circle in $\Sigma$. The polarity with respect to $\Sigma$ maps the plane of $\bar{k}$ to a straight line $k^*$. The points $S^*$ of $k^*$ determine M-spheres $S$ through $k$. These spheres are said to form a {\em pencil of %M-
 spheres}.

\begin{figure}[htb]
\centering{%\vspace{4cm}
\includegraphics%[width=0.9\textwidth]
{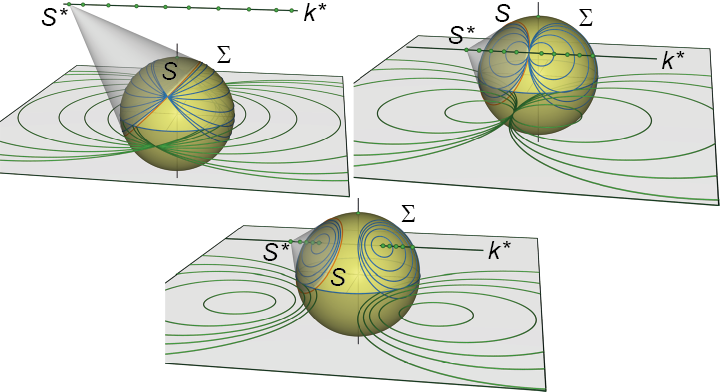}
}
 \caption{Pencils of %M-
 circles in 2D M\"obius geometry are stereographic projections of pencils of circles $S$
 on $\Sigma$; the latter are easy to understand since they are cut out from $\Sigma$ by a pencil of planes.
 The poles $S^*$ of these planes lie on a straight line $k^*$. The figure illustrates the $3$ types of pencils, which by rotation about the central line generate the $3$ types of pencils of %M-
 spheres. }
	\label{fig:pencils}
\end{figure}

 A \emph{pencil} is more generally defined as the set of all M-spheres $S$ whose corresponding points $S^*$ form a straight line $k^*$ in $P^4$.
 If $k^*$ does not intersect $\Sigma$, we obtain an {\it elliptic} pencil, formed by all M-spheres through a real circle. If $k^*$ is tangent to $\Sigma$, we get a {\it parabolic} pencil, and if $k^*$ intersects $\Sigma$ in two %real
 points, we have a {\it hyperbolic} pencil. A parabolic pencil is formed by all M-spheres touching a plane $P$ at a fixed point $\pw \in P$. The common set of a hyperbolic pencil is an imaginary circle; the pencil contains
 two null spheres, %of radius zero,
 corresponding to the points $\Sigma \cap k^*$. Clearly, in all cases the centers of the M-spheres lie on a
 straight line $A$, which is also the image of the line $k^*$ under the extended stereographic projection.
 Due to the rotational symmetry of pencils, it is sufficient to understand the planar sections
 through the rotational axis $A$, which are pencils of %M-
 circles (see Fig.~\ref{fig:pencils}).

 {\em Linear families of %M-
 spheres} $S$ are those where the set of poles $S^*$ is a $d$-dimensional subspace  of $P^4$.
 Apart from the already discussed pencils ($d=1$), we still have to consider the cases $d=2$ ({\em bundles of %M-
 spheres}) and $d=3$ ({\em
 %M-
 sphere complexes}).

 Before doing so, we should look at {\em orthogonal M-spheres}, i.e., those which intersect under a right angle. Two real spheres $S_1, S_2$ with centers $\mw_1, \mw_2$ and radii $r_1, r_2$, are orthogonal if
 $(\mw_1 - \mw_2)^2=r_1^2+r_2^2$. Their corresponding points $S_i^*$ (in pentaspherical coordinates) are
 $$ S_i^*=(2m_{i,x}, 2m_{i,y}, 2m_{i,z}, \mw_i^2-r_i^2-1, \mw_i^2-r_i^2+1).$$
 Two points $A=(a_1,\ldots,a_5)$ and $B=(b_1,\ldots,b_5)$ in $P^4$ are conjugate with respect to $\Sigma$ (one lies in the polar hyperplane of the other one) if $\langle A,B \rangle:= a_1b_1+ \ldots + a_4b_4-a_5b_5=0$. We see immediately that $\langle S_1^*, S_2^* \rangle = 0$ if and only if $(\mw_1 - \mw_2)^2=r_1^2+r_2^2$. Hence, {\it orthogonal real spheres $S_1, S_2$ are characterized by conjugate points} $S_1^*, S_2^*$. The latter characterization allows to define \emph{orthogonality} also for null and imaginary spheres.

\begin{figure}[htb]
\centering{%\vspace{4cm}
\includegraphics[width=0.5\textwidth]{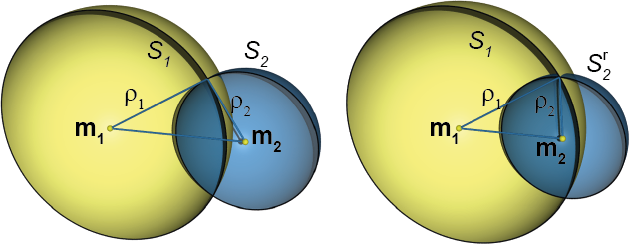}
}
 \caption{Orthogonality between M-spheres $S_1, S_2$ is not only defined for real spheres (left), but also for imaginary ones. A real sphere $S_1$ is orthogonal to an imaginary sphere $S_2$ (with center $\mw_2$, radius $i\rho_2$), if it intersects the real representer $S_2^r$ (with center $\mw_2$, radius $\rho_2$) of $S_2$ in a great circle (right).}
	\label{fig:orthogonal-spheres}
\end{figure}

 \begin{rem} \label{rem:complex-sphere}
 In the above discussion, it can happen that one of the spheres, say $S_2$, is imaginary,
 %has a real center, but an imaginary radius $r_2=i \rho_2$; such complex spheres
 and thus corresponds to a point $S_2^* \in \R^4$  in the interior  of the $3$-sphere $\Sigma$. Hence, the orthogonality condition now reads $(\mw_1 - \mw_2)^2=r_1^2-\rho_2^2,$ where $\rho_2^2:=-r_2^2$.
 %Let us define the {\em real representer} $S_2^r$ of the \mscomm{complex} sphere $S_2$ as the sphere with center $\mw_2$ and radius $\rho_2$.
 Then, the two spheres $S_1$ and $S_2$ are orthogonal if and only if the intersection circle of $S_1$ and $S_2^r$ is a great circle of $S_2^r$ (see Fig.~\ref{fig:orthogonal-spheres}).
 \end{rem}

It is now easy to understand a {\em linear complex of %M-
spheres}.  The points $S^*$ corresponding to the M-spheres of the complex lie in a hyperplane $H$, whose pole shall be denoted by $A^*$.  The point $A^*$ is conjugate to all points $S^* \subset H$
and it defines a (not necessarily real) M-sphere $A$, which is orthogonal to all M-spheres $S$ of the complex (and which is seen in the spherical model as intersection $\bar{A}=H \cap \Sigma$). Hence, {\it a complex is formed by all M-spheres orthogonal to a fixed M-sphere $A$}. %(see Fig.~\ref{fig:sphere-complexes}).
We have again $3$ types, depending on whether $A$ is real, imaginary or null, i.e.,
 degenerates to a point.
In the latter case, all M-spheres of the complex pass through this point.

%\begin{figure}[htb]
%\centering{\bf Figure to be removed
%\includegraphics[width=9.624cm]{Figure5.png}
%}
%\caption{The three types of sphere complexes arise as sets of spheres which are (i) orthogonal to a real sphere (left), (ii) orthogonal to a complex sphere (middle) or (iii) pass through a fixed point (right).}
%\label{fig:sphere-complexes}
%\end{figure}

\begin{figure}[htb]
\centering{%
\includegraphics[width=.6\textwidth]{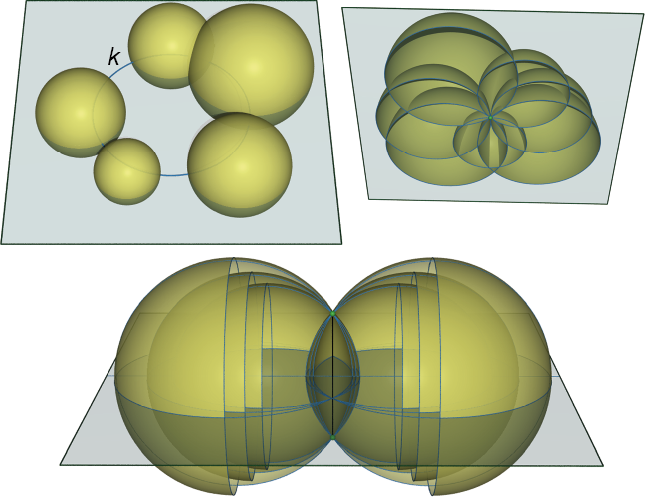}
}
 \caption{The three types of %M-
 sphere bundles are defined as sets of M-spheres which intersect an M-circle $k$ under right angle. Top left: for an elliptic bundle, $k$ is real. Bottom: in a
 hyperbolic bundle, $k$ is imaginary, implying that the spheres pass through $2$ fixed points. Top right: a parabolic bundle is a limit case and formed by all M-spheres centered in a plane and passing through a fixed point in that plane.}
 \label{fig:sphere-bundles}
\end{figure}

A {\em bundle of %M-
spheres} is defined by a 2D plane $P \subset P^4$, which may be viewed as the intersection of two hyperplanes $H_1, H_2$. Hence, the M-spheres of a bundle belong to two complexes
and thus they are orthogonal to two (not necessarily real) M-spheres $A_1, A_2$. This requires orthogonality to the (not necessarily real) intersection M-circle $A_1 \cap A_2$,
and thus the M-spheres of a bundle have their centers in a plane and they are all orthogonal to an M-circle $k$ in that plane (note that the M-circle $k$ is  seen in the spherical model as the intersection M-circle $\bar{k}= P \cap \Sigma $). Once again there are three types, since $k$ may be real, imaginary or null
%degenerate to a point
(see Fig.~\ref{fig:sphere-bundles}). The case of a imaginary circle is easy to understand,
since now  all  M-spheres of the bundle pass through two real points (on the real axis of the imaginary circle). This follows for example from the fact that in this case two hyperplanes through $P$ are tangent to $\Sigma$ (or by observing that
the polar dual to the plane $P$ is a straight line, determining a hyperbolic pencil of %M-
spheres which are all orthogonal to the %M-
spheres of the bundle).

 {\em M\"obius transformations} (or \emph{M-transformations}) act bijectively on the points and spheres of M\"obius geometry. M-transformations are the only conformal (angle preserving) maps in Euclidean $3$-space; they include Euclidean similarities as special cases. In the spherical model they are seen as those projective transformations $\kappa\colon P^4\to P^4$ which map $\Sigma$ onto itself. Using homogeneous coordinates, the equation of $\Sigma$ is $ \Xw^T \cdot J \cdot \Xw =0$, with the diagonal matrix $J={\rm diag}(1,1,1,1,-1)$. A projective map $\kappa:\ \Xw \mapsto A \cdot \Xw$ in $4$-space maps $\Sigma$ onto itself  and thus induces an M-transformation if and only if its matrix $A$ satisfies
\be \label{M-trafo}
   A^T\cdot J \cdot A = \lambda J,
\ee
for some $\lambda \ne 0$.

The basic example of an M-transformation is the {\it inversion with respect to an M-sphere} $S$.  In the spherical model it appears as a projective symmetry of $\Sigma$, more precisely, a harmonic perspective collineation $\kappa$ which maps $\Sigma$ onto itself. The axis of this
collineation is the hyperplane $H_s \supset \Sb$ and the center is its pole $S^*$. Thus %We can also define
an inversion with respect to an imaginary sphere $S$ is also defined. In this case, $S^*$ lies in the interior of $\Sigma$. In $3$-space,
such an inversion can be realized by reflection at the (real) center of $S$ and subsequent inversion with respect to the real representer $S^r$ of $S$. %(see Remark \ref{rem:complex-sphere}).
Notice that the real representer is not \emph{M-invariant}, i.e., $\kappa (S^r)\ne (\kappa S)^r$ in general.

%% Commented out by Mikhail:
%{\it Maybe we can illustrate the power of M-transformations to bring the pencils, bundles and complexes from above into special Euclidean forms. For pencils this would be: elliptic: pencil of lines; parabolic: pencil of parallel lines; hyperbolic: pencil of concentric circles. This may be useful later for design as well.}

 \subsection{Darboux cyclides in the spherical model of M\"obius geometry.} \label{ssec:pencil}
 %---------------------------------------------------------------------------
 %
In this subsection we will elaborate on the classical result that {\em a Darboux cyclide $D$ is the stereographic projection of a surface $\bar{D}$ in $\Sigma$ which is obtained as the intersection of $\Sigma$ with another quadric $\Gamma$} (see \cite{coolidge:1916}; here by definition the ideal point $\infty$ belongs to a nontrivial cyclide $D$ if and only if $\lambda=0$ in equation~\eqref{D}). In other words,  $\bar{D}$ is the common set of points, the so-called {\it carrier}, of a pencil of quadrics which contains the sphere $\Sigma$. Deriving this result is not hard, and it will give us important insight needed later in our studies.

Let us first show that inverse stereographic projection maps a Darboux cyclide $D$ to the carrier of a pencil of quadrics. For that, we substitute from (\ref{stereo2}) into (\ref{D}). Splitting  $Q$
%and $L$
into homogeneous terms of degree
$2$, $1$, and $0$ as $Q(\xw)=q(\xw)+\qw \cdot \xw + q_0$ and writing $L=\lw\cdot \xw$
%$L=\lw\cdot \xw + l_0$
yields with $\bar{x}:=(x_1,x_2,x_3)$
$$ \lambda\left({\bar{x}^2 \over (1-x_4)^2}\right)^2+  {\bar{x}^2 \over (1-x_4)^2}%(
{\lw \cdot \bar{x} \over 1-x_4}%+l_0)
+ {q(\bar{x}) \over (1-x_4)^2}
+{\qw \cdot \bar{x} \over 1-x_4}+q_0=0.$$
Inserting the equation of $\Sigma$, $\bar{x}^2=1-x_4^2$, we finally obtain the equation of a quadric or a plane in $4$-space,
$$ \Gamma: \ \lambda(1+x_4)^2 + (1+x_4)%(
\lw \cdot \bar{x} %+ l_0(1-x_4))
+ q(\bar{x})+\qw \cdot \bar{x}(1-x_4) + q_0(1-x_4)^2=0. $$
We may assume without loss of generality that $\Gamma$ is a quadric by adding a multiple of $\bar{x}^2+x_4^2-1$ to the left-hand side of the equation.
The quadric $\Gamma$ contains the projection center $Z$ if and only if $\lambda=0$, i.e., $\infty\in D$.

Conversely, we have to show that intersection of $\Sigma$ with any quadric $\Gamma \subset \R^4$ yields a surface $\bar{D}$ whose
stereographic image in $3$-space is of the form (\ref{D}), i.e., a Darboux cyclide. Again separating terms of different degrees and splitting
$(x_1,\ldots,x_4)$ as $(\bar{x}, x_4)$, we write the equation of the quadric as
$$ \Gamma: q_2(\bar{x}) + q_1(\bar{x})x_4 + q_0 x_4^2 + l_1(\bar{x})+l_0x_4 + c_0 =0,$$
where $q_2$, $\{ q_1,l_1 \}$ and $\{q_0,l_0,c_0\}$ are of degrees $2$, $1$, and $0$, respectively.
Insertion from~(\ref{stereo1}) yields
$$ (q_0+l_0+c_0)(\xw^2)^2+ 2 \xw^2[q_1(\xw)+l_1(\xw)]+4q_2(\xw)+Q(\xw)=0,$$
with
$$ Q(\xw)=2(c_0-q_0)\xw^2+2(l_1(\xw)-q_1(\xw))+q_0-l_0+c_0. $$
This equation is of the form (\ref{D}) and thus represents a Darboux cyclide.

Transferring the study of Darboux cyclides to the study of the carriers of pencils of quadrics which
contain the sphere $\Sigma$ has a lot of advantages since we are now able to apply knowledge from projective
geometry on pencils of quadrics.
%For example, this reduces the classification of cyclides to the classification of pencils of quadrics in $4$-space.
The study of all possible types of cyclides, i.e., pencils of quadrics in $4$-space, can be elegantly based on the Weierstrass theory of elementary divisors
(a workaround being found in \cite{takeuchi:2000}; see also \cite[p.~301]{coolidge:1916} and \cite{uhlig:1976}).

%\mscomm{The conformal classification of cyclides is as follows.}

%\mscomm{\textbf{Added the following theorem for completeness and also because the Takeuchi statement is %not true literally for our (more general) definition of cyclides. M.S. }}

%\begin{thm} \cite[Section~0 and Theorem~6.6]{takeuchi:2000} Any cyclide can be brought by a (real) %M\"obius transformation to one of the surfaces:
%\begin{align*}
%&\lambda(x^2+y^2+z^2)^2+\mu x^2+\nu y^2+\kappa z^2+\eta=0
%   \qquad\text{or}\\
%&\mu x^2+\nu y^2+\kappa z=0.
%\end{align*}
%\end{thm}

%\mscomm{Mikhail: I've obtained a short proof of the Takeuchi conformal classification modulo Uhlig result --- I can %add this if necessary.}

\subsection{Families of circles on Darboux cyclides.} \label{ssec:circles}
%---------------------------------------------------
%
Working in the spherical model, we represent the Darboux cyclide  as the carrier $\Db$ of the
pencil $\Pz$ of quadrics
\be \Xw^T \cdot (A-tJ) \cdot \Xw =0, \label{pencil} \ee
spanned by the sphere $\Sigma: \Xw^T \cdot J \cdot \Xw =0$, and another quadric $\Gamma:\Xw^T \cdot A \cdot \Xw =0$.
We can work with the inhomogeneous pencil parameter $t$, since we do not have to represent $\Sigma$.
%; the minus sign is for later convenience.
Obviously, given a point $Y \notin \Sigma$ (with coordinates $\Yw$), there is
a unique quadric $\Gamma_Y$ in $\Pz$ which contains $Y$, namely the one with the pencil parameter
$t=(\Yw^T \cdot A \cdot \Yw) / (\Yw^T \cdot J \cdot \Yw)$.

If the carrier $\Db$ of $\Pz$ contains a circle $\kb$ in a plane $P_k$, we can choose
any point $Y \notin k$ in $P_k$ and investigate the quadric $\Gamma_Y \in \Pz$ which
contains $Y$. The intersection of the plane $P_k$ with $\Gamma_Y$ must contain the circle $\kb$ and
the point $Y$ and thus the entire plane $P_k$ must be contained in the quadric $\Gamma_Y$. A regular
quadric in 4-space can contain at most straight lines, but no planes. Hence, $\Gamma_Y$ is
singular, i.e., a {\em quadratic cone}.

Hence, looking for circles of $\Db$, we have to search for their planes on the quadratic cones in $\Pz$.
This is a well-known subject, and can be analytically accessed as follows. A quadric in the
pencil for some parameter $t$ is a cone if $$\det(A-tJ)=0.$$ Hence, the possible parameter
values $t_i$ are the \emph{real} eigenvalues of the matrix $A \cdot J^{-1} = A \cdot J$. Let $t_i$ be
such an eigenvalue. Then the linear system $(A-t_iJ) \cdot \Xw=0$ has a nontrivial solution, say
$\Vw_i$, which describes a point $V_i$, the vertex of the cone. This implies
$$ A \cdot \Vw_i = t_i J \cdot \Vw_i,$$
which says that the polar hyperplanes of the point $V_i$ with respect to the quadric
$\Gamma \in \Pz: \Xw^t \cdot A \cdot \Xw =0$
and with respect to the sphere $\Sigma$ agree. In fact, the polar hyperplanes of $V_i$ with respect to
all quadrics in $\Pz$ are identical, because of
$$ (A-tJ)\cdot \Vw_i = (t_i-t)J\cdot \Vw_i.$$
We see that $V_i$ is a vertex of the common polar simplex of the pencil $\Pz$ (which needs
not be a regular simplex, needs not have only real vertices and may degenerate in various ways).
In the generic case, we get
$5$ linearly independent eigenvectors of the matrix $A\cdot J$, determining the $5$ vertices of
the common polar simplex. Only real vertices $V_i$ are of interest to us. Connecting any of these vertices $V_i$ with the carrier $\Db$, we
obtain a quadratic cone $\Gamma_i$ which contains $\Db$.

\begin{figure}[htb]
\centering{%\vspace{4cm}
\includegraphics[width=.6\textwidth]{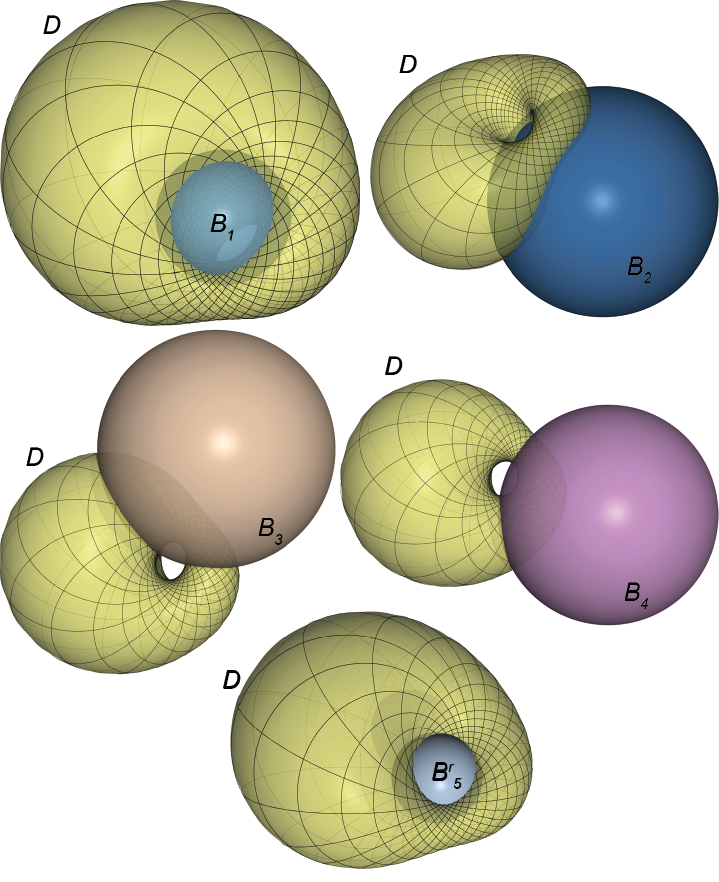}
}
 \caption{A cyclide $D$ determines a pencil of quadrics in $4$-space. Each vertex $V_i$ of the common polar simplex of the pencil corresponds to a sphere $B_i$ in $3$-space with respect to which $D$ is symmetric (inversion
 with respect to $B_i$ maps $D$ onto itself). This figure shows a cyclide $D$ which is symmetric with respect
 to $5$ (pairwise orthogonal) spheres $B_1,\ldots, B_5$; since $B_5$ is imaginary, we illustrate
 its real representer~$B_5^r$.}
 %\mscomm{(inversion in $B_5$ is by definition inversion in $B_5^r$ composed with the central symmetry with respect to the midpoint of $B_5^r$).}}
 \label{fig:M-symmetry}
\end{figure}

\begin{rem} \label{M-symmetry} Any real vertex $V_i \notin \Sigma$ of the common polar simplex, together with its
polar hyperplane, define a joint projective symmetry of all quadrics in the pencil, and thus of its carrier $\Db$.
The vertex $V_i$ determines an M-sphere, say $B_i$, in $3$-space (in the notation used above,
we have $B_i^*=V_i$) . The inversion with
respect to such a sphere $B_i$ maps the cyclide $D$ onto itself. Since the vertices $V_i$ are pairwise conjugate with respect to $\Sigma$, the
spheres $B_i$ are pairwise orthogonal. In the case with $5$ real vertices $V_i$, we get cyclides $D$ which are M-symmetric with respect to $5$ pairwise orthogonal M-spheres (see Fig.~\ref{fig:M-symmetry} and \cite{coolidge:1916});
note that one of the $5$ points $V_i$ must lie in the interior of $\Sigma$ and thus the corresponding
sphere is imaginary ($B_5$ in Fig.~\ref{fig:M-symmetry}).
\end{rem}

\begin{figure}[htb]
\centering{
\includegraphics[width=.4\textwidth]{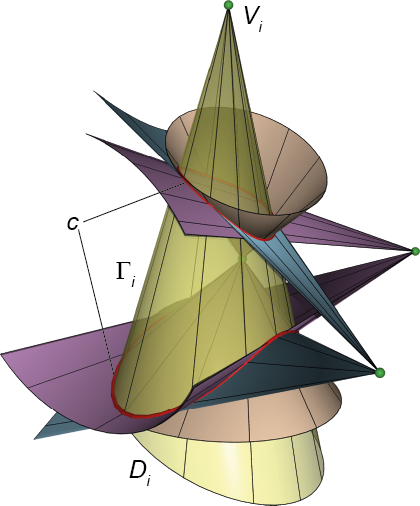}
}
 \caption{A pencil of quadrics in $3$-space has an algebraic curve $c$ of order $4$ as its carrier.
 There are four cones $\Gamma_i$ in this pencil; projecting $c$ from a cone vertex $V_i$ onto a plane $H_i$ results in a conic $D_i$.}
 \label{fig:quadpencil}
\end{figure}

Once again we look at a lower-dimensional counterpart to visualize the situation. Fig.~\ref{fig:quadpencil}
shows a pencil of quadrics in $3$-space, spanned by a sphere and another quadric. The pencil contains four
cones with the vertices at the vertices of the common polar simplex. The figure also shows that the projection of the intersection curve
(algebraic curve $c$ of order four) from a point $V_i$ onto a plane is a conic $D_i$ and that not all points of $D_i$
need to be projections of real points of $c$.

\begin{figure}[htb]
\centering{
\includegraphics[width=.65\textwidth]{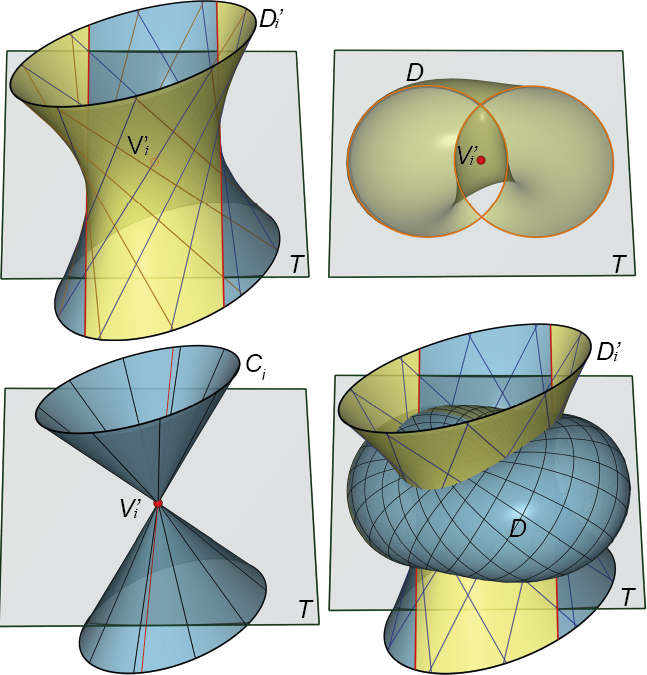}}
 \caption{Extended stereographic projection of the situation in 4D. (Top left) Connecting the image point $V_i'$ of the cone vertex $V_i$ with the rulings of the projected ruled quadric $D_i'$ results in planes $T$ which carry circles of the cyclide $D$. (Top right) Each %cone tangent
 plane $T$ %obtained in this way
 carries two circles and is tangent to the cyclide $D$ at $\le 2$ points, the common points of these two circles. (Bottom left)
 %The circles of a paired family lie in the tangent planes
 Each plane $T$ is tangent to a quadratic cone $C_i$ with vertex $V_i'$ which is tangent to $D_i'$. (Bottom right) A paired family of circles arising from the two families of rulings in $D_i'$.
 }
 \label{quadprojection}
\end{figure}

Let us return to $4$-space. We may project $\Db$ from a cone vertex $V_i$ onto a hyperplane $H_i$ (not passing
through $V_i$) and obtain a quadric
$D_i \subset H_i$.
The structure of this projection becomes clear in $\mathbb{C}^4$. %, we switch to this space for a moment.}
The quadric $D_i$ appears as ``doubly covered'', i.e., each of its points is the image of two points of $\Db$. However, not any real point of the quadric
$D_i$ is the image of two real points on $\Db$ (see Fig.~\ref{fig:quadpencil}); it can as well be the image of two conjugate complex points of $\Db$. The double coverage of $D_i$ explains why it can be the image of a $2$-surface
$\Db$ of algebraic order $4$ under a central projection.
Returning to $\mathbb{R}^4$ we get that
if the cone $\Gamma_i$ contains planes, they pass through circles on $\Db$ and of course through $V_i$ and thus
these planes (and the circles) are projected onto straight lines of the quadric $D_i$. This
shows us how to obtain circles on Darboux cyclides (see Fig.~\ref{quadprojection}). We summarize our findings in the following theorem which is a classical result, but probably not so straightforward to extract from the classical texts (e.g., \cite{coolidge:1916}) so that it did not get used in later studies on cyclides (\cite{blum:1980,takeuchi:2000}).

 \begin{thm} In the spherical model of M\"obius geometry, a Darboux cyclide $D$ appears as the
 carrier $\Db$ of a pencil $\Pz$ of quadrics which contains the sphere $\Sigma$. Any real vertex of the
 (possibly degenerate) common polar simplex of the quadrics in $\Pz$ is the vertex of a quadratic
 cone through $\Db$. If such a cone contains planes, these planes intersect $\Sigma$ in circles of $\Db$, which
 determine circles or lines on $D$.
 Any circle or line on $D$ can be found in this way.
 \end{thm}

This result gives rise to several conclusions on circle families on nontrivial irreducible Darboux cyclides (containing no spheres or planes as a whole).
%(but the converse is not true in general. %\mscomm{Hereafter we assume that the cyclides are \emph{irreducible}, i.e., contain no spheres or planes as a whole, and \emph{not completely singular}, i.e., contain at least one nonsingular real point}.
%splitting into a union of spheres or planes.
There are two types of quadratic cones in the $4$-space which contain planes (but no hyperplane):

\begin{enumerate}
\item The cone $\Gamma_i$ is obtained
by connecting a nonsingular ruled quadric $D_i$ in a hyperplane $H_i$ with a point $V_i$ outside $H_i$. The
ruled quadric $D_i$ carries two families of straight lines (\emph{rulings}). Each plane through a ruling of one
family intersects $D_i$ in a ruling of the second family. Connecting with the vertex $V_i$ and
interpreting the scene from the M\"obius geometric perspective, we find two circle families on $\Db$ (and $D$).
We call them
{\em paired families}, since one family determines the other one. Any sphere through a circle of one family
intersects the cyclide $D$ in a second circle which belongs to the other family. This holds in particular
for the plane of such a circle: it intersects the cyclide in a circle of the other family (see Fig.~\ref{quadprojection}). We shall say that paired families are \emph{special}, if $V_i\in\Sigma$ --- this is only possible if $D$ is the image of a quadric under an inversion
\cite[cf.~Theorem~35 in p.~301]{coolidge:1916}. Paired families were called \emph{residual} in \cite{coolidge:1916}.

\begin{figure}[htb]
\centering{
\includegraphics[width=.5\textwidth]{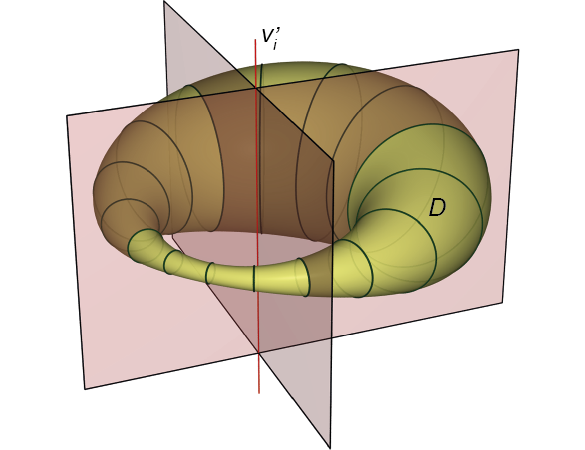}}
 \caption{Single families of circles are found only on those Darboux cyclides which are envelopes of a one-parameter family
 of spheres. The spheres touch the cyclide $D$ along the circles of the single family, and the planes of the circles pass through a fixed
 line $v_i'$ (image of the vertex line $v_i$ of a quadratic cone containing $\Db$).}
 \label{quadprojection2}
\end{figure}

\item The cone $\Gamma_i$ is the connection of an irreducible quadratic cone $D_i$ (with vertex $W_i$) in a hyperplane $H_i$ and a point $V_i$ outside $H_i$. This is the same as connecting a conic (on $D_i$) in some plane $P \subset H_i$ with a straight line $v_i=V_iW_i$ which does not intersect $P$. The straight line is the vertex line of the cone. Such a cone %\mscomm{(if its real points does not degenerate to a single plane)}
    carries a single family of planes (through $v_i$) and thus we speak of a {\em single family of circles} on $\Db$ and on the underlying Darboux cyclide $D$. Obviously, a sphere through a circle of a single family intersects the cyclide in a circle of the same family, just like a plane through a ruling of the quadratic cone $D_i$ intersects the cone in another ruling. However, there are also tangent planes of the cone $D_i$. They lead to spheres touching the cyclide along the circles of a single family and thus we see that a Darboux cyclide which contains a single family is a {\em canal surface}, i.e., the envelope of a one-parameter family of spheres (see Fig.~\ref{quadprojection2}). In fact, almost each such Darboux cyclide is the image of a \emph{complex} quadratic cone under a \emph{complex} inversion \cite[Theorem~38 in p.~302]{coolidge:1916}.
    \end{enumerate}

 \begin{prop} \label{prop-families} There are two different types of circle families on nontrivial irreducible Darboux cyclides $D$. (i) Paired families are two families
 such that a sphere through a circle of one family intersects the cyclide in another circle, which belongs to the second
 family. The planes of both circle families are either the tangent planes of a quadratic cone or the planes of two
 pencils with intersecting axes (in the sense of projective geometry). (ii) The circles of a single family are the characteristic circles of a generation of the cyclide as a canal surface; a sphere through a circle of a single family either touches the cyclide along that circle or it contains another circle of $D$ belonging to the same family. The planes of the circles lie in a pencil.
 \end{prop}

 The statements on the planes of circles follow from the extended stereographic projection of the plane family in the quadratic cone
 $\Gamma_i$ (see Fig.~\ref{quadprojection}). Indeed, we may assume that $Z\ne V_i$ and $Z\in H_i$. If $Z\not\in D_i$ in type (i) then the lines in the irreducible quadric $D_i$ are projected to tangents to a certain conic in $3$-space. Joining the resulting lines with the projection of the point $V_i$ we get tangent planes of a quadratic cone. The case of two pencils with intersecting axes in type~(i) arises if the projection center $Z$, the vertex $V_i$, and a point of the ruled base quadric $D_i$ are collinear. This is only possible if $\lambda=0$ in \eqref{D};
  an example is the paired family of circles in an ellipsoid (not of revolution).
%{\bf We should make an example for this and illustrate.}

 %\mscommnew{\bf Needed for the sequel:}

 Using this approach one can obtain further information about the arrangement of circle families on Darboux cyclides. Once the quadric $D_i$ contains at least one line, it contains a line through each point. Thus any family of circles on an irreducible cyclide $D$ covers the entire cyclide. Now assume that there are two circle families on $D$. Then through a generic point $p\in D$ one can draw circles from both families. There are two possibilities. First, these circles may have another common point (or touch each other). In this case they lie in one sphere. Hence by Proposition~\ref{prop-families} the circles belong to two paired families. Second, the point $p$ may be the unique transversal intersection point of the circles. In this case any two other circles from the families have also a single intersection point because the parity of the number of intersection points (counted with multiplicities) does not depend on the choice of two particular circles from the families.

 \begin{prop} \label{prop-1commonpoint} In nontrivial irreducible Darboux cyclides, two circles from distinct families, which are not paired together, intersect each other transversally at a single point.
 \end{prop}

 \begin{rem} Similarly, there are two types of families of \emph{conics} in the intersection $\bar D$ of two arbitrary nonsingular $3$-dimensional quadrics in $4$-space ($\bar D$ containing no $2$-dimensional quadrics as a whole). \emph{Paired} families are such that a hyperplane through a conic of one family intersects $\bar D$ at another conic, which belongs to second family. A hyperplane through a conic of a \emph{single} family either touches $\bar D$ or intersects it in another conic of the same family. This is proved analogously to Proposition~\ref{prop-families}.
 %Propositions~\ref{prop-1commonpoint}, \ref{prop-families} \mscomm{(without the assertions on the planes of the circles and generation of the surface as a canal surface)}, and their proofs remain true, if one replaces ``\mscomm{nontrivial irreducible} Darboux cyclides'' by ``intersections of two nonsingular \mscomm{$3$-dimensional} quadrics in $4$-space \mscomm{containing no $2$-dimensional quadrics as a whole}'', ``circle'' by ``conic'', and \mscomm{``sphere'' by ``hyperplane''}.
 %In particular, this result allows to describe the families of isotropic circles in so-called isotropic cyclides.
 \end{rem}

 \subsection{The number of families of real circles on a Darboux cyclide.} \label{ssec:realcircles}
 %-------------------------------------------------------------------------
 Let us further discuss the circle
 families on cyclides and start with some simple examples. Darboux cyclides which contain two single families must be the envelopes of two families of spheres and thus these are the {\em Dupin cyclides} which have already
 found a lot of attention in Geometric Modeling (see, e.g., \cite{degen:2002,krasauskas-maurer:2000}). As a special case, consider a ring torus, which
 obviously contains two single families (meridian circles and the circular paths of the generating
 rotation); furthermore it carries a paired family, namely the two families of Villarceau circles. Likewise any Dupin ring cyclide carries four families of circles, since it is the image of a ring torus under an inversion.
   Rotating a circle about an axis (not meeting the axis of the circle and not contained in the plane of the circle) we may obtain a ring shape, which is a
 Darboux cyclide carrying five real families of circles (two paired ones and the single family of rotational circles;
 for in illustration, see
 %\mscomm{Fig.~\ref{fig:web-type4} and}
 \cite{archgeom:2007}).

\begin{figure}[htb]
\centering{
\includegraphics[width=.6\textwidth]{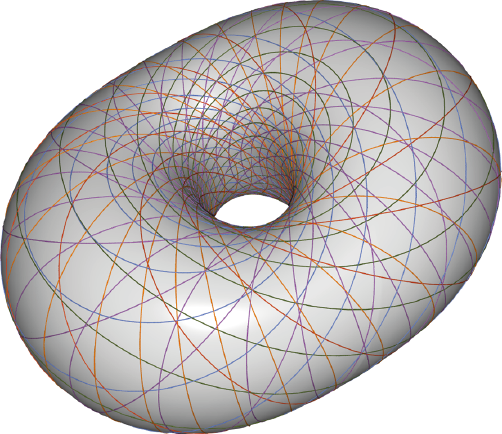}}
 \caption{A cyclide with $6$ families of circles, shown in different colors. We will see later that we can form
 eight different $3$-webs from these circles (cf. Figure \ref{fig:closed-web}).}
 \label{fig:6families}
\end{figure}

 However, there are also {\em Darboux cyclides which carry  six families of real circles}. This has first been pointed out by \cite{blum:1980}, at hand of a special case of symmetric cyclides, and in a quite involved elementary proof which gives
 less insight than our approach. Blum mentions that Darboux had pointed to 10 families of circles, but did not differentiate
 between real and complex ones.  \cite{takeuchi:2000} provides a discussion of the real families of circles
 on cyclides, similar to our approach, but he does not address the geometric relation to the cones in the pencil of quadrics.
 In none of these papers, we find a figure of such a surface, which of course is easy to generate with our approach and
 provided, e.g., in Fig.~\ref{fig:6families}.

We can easily obtain these results and more as follows.
Without caring about reality and looking only at the generic case, we get Darboux's $10$ families, two families from each of the $5$ quadratic cones.  To determine the number of real circle families on a cyclide we
have to investigate the quadratic cones in a pencil $\Pz$ of quadrics through $\Sigma$. Using the notation from equation~(\ref{pencil}), the matrices of the quadratic cones $\Gamma_i$ are $C_i:=A-t_iJ$, where $t_i$ is a real eigenvalue of the matrix $A\cdot J$.
A quadratic cone with two families of planes (which determines a paired family of circles) is given by a matrix $C_i$ of signature $(0,+,+,-,-)$. For a cone determining a single family, the matrix has signature either $(0,0,+,+,-)$ or $(0,0,+,-,-)$.

For brevity, let us now focus on the generic case with $5$ distinct vertices of the common polar simplex, where we also assume that all five vertices are real. In this case, we can apply a projective map in $4$-space which maps these vertices to the origin and the ideal points of the four coordinate axes. In other words, we have achieved that matrix $A$ is a diagonal matrix, say with diagonal $(a_1,\ldots,a_5)$. Matrix $A \cdot J$ is then also diagonal with eigenvalues $t_i=\{a_1,\ldots,a_4,-a_5\}$, and so the eigenvalues of $C_i$ are $(a_1-t_i,\ldots, a_4-t_i,a_5+t_i)$.
%It is very easy to see that the smallest and largest of the $t_i$'s  will not generate a  sequence $(0,+,+,-,-)$, hence at most $3$ cones can lead to paired families. Thus we will never have more than $6$ families of real circles on a Darboux cyclide of the considered type.

\begin{exmp}
To give an example of a cyclide with $6$ real circle families (Fig.~\ref{fig:6families}), we consider %define
the pencil with the matrix $A={\rm diag}(-2,-1,1,2,0)$, which defines a quadratic cone  $\Gamma_1$ with vertex at the origin; $\Gamma_1$ contains two plane families. Using Cartesian coordinates, our pencil is spanned by
  $$ x_1^2+\ldots+x_4^2-1=0, \ -2x_1^2-x_2^2+x_3^2+2x_4^2=0.$$
 Following the above procedure, or equivalently, eliminating one of the variables, we arrive at the following $4$ additional
 cones (cylinders parallel to the coordinate axes):
 $$ \Gamma_2 (t_i=-2):\ x_2^2+3x_3^2+4x_4^2-2=0, \ \Gamma_3 (t_i=-1):\ -x_1^2+2x_3^2+3x_4^2-1=0,$$
 $$ \Gamma_4 (t_i=1):\ -3x_1^2-2x_2^2+x_4^2+1=0, \ \Gamma_5 (t_i=2):\ -4x_1^2-3x_2^2-x_3^2+2=0.$$
 We see that $\Gamma_3$ and $\Gamma_4$ determine paired families, while the largest and smallest $t_i$ determine cylinders with an oval base quadric and thus do not contribute further real circles.
 \end{exmp}

\begin{rem} \label{rem:largest eigenvalue}
For the smallest real eigenvalue $t_0$ of the matrix $A \cdot J$ the matrix $A-t_0J$ cannot have signature $(0,+,+,-,-)$. Indeed, the eigenvalues of the matrix $A-tJ$ change continuously with decreasing $t$, but for $t<t_0$ they should have the same signs $(+,+,+,+,-)$ as those of~$J$. Similarly, the largest real eigenvalue does not lead to a paired family. Thus there are at most $3$ paired families. Single families may only appear for eigenvalues of multiplicity $2$, hence we never have more than $6$ families of real circles on a nontrivial irreducible Darboux cyclide.
\end{rem}

This can easily be cast into a simple {\em algorithm for computing the circle families on a given cyclide}.
Assuming that we have already identified the corresponding pencil of quadrics~(\ref{pencil}), we
%and computed a real eigenvalue $t_i$ of the matrix $A \cdot J$.
compute the real eigenvalues of the matrix $A \cdot J$.
For each such eigenvalue $t_i$, we compute the eigenvalues of the matrix $A-t_iJ$. If the eigenvalues of the matrix $A-t_iJ$ have signs either $( +, +, -, -, 0)$ or $(+,+,-,0,0)$ or $(+,-,-,0,0)$ then
%denote them by $(k_1^2,  k_2^2,  -k_3^2,  -k_4^2,  0) $ with $k_1, \ldots, k_4 > 0 $.
%The square denotation is for latter convenience.
we bring the matrix to the form $\mathrm{diag}(k_1^2,k_2^2,-k_3^2,-k_4^2,0)$
by a linear transformation $\kappa$ of coordinates,
where $k_1,k_3>0$ and $k_2,k_4\ge 0$. Denote by $P_1(s)$ and $P_2(s)$, where $s$ is a parameter,
the images of the planes
$$
\begin{cases}
k_1x_1+sk_2x_2-k_3x_3-sk_4x_4=0,\\
sk_1x_1-k_2x_2+sk_3x_3-k_4x_4=0;
\end{cases}
%\qquad
\text{ and }
%\qquad
\begin{cases}
k_1x_1-sk_2x_2-k_3x_3-sk_4x_4=0,\\
sk_1x_1+k_2x_2+sk_3x_3-k_4x_4=0;
\end{cases}
$$
under the transformation $\kappa^{-1}$. In particular,
the families $P_1(s)$ and $P_2(s)$ coincide, if either
$k_2=0$ or $k_4=0$. Stereographic projections of $P_1(s)\cap\Sigma$ and $P_2(s)\cap\Sigma$ are families of circles on the cyclide.
If the eigenvalues of $A-t_iJ$ have signs neither $( +, +, -, -, 0)$ nor $(+,+,-,0,0)$ nor $(+,-,-,0,0)$
%in certain order
then $t_i$ does not contribute to circle families.

 \subsection{The polar dual.}  \label{ssec:polardual}
 %-----------------------------
 %
 While not absolutely necessary to get further insight, we think that it is helpful
 to also consider the polarity with respect to $\Sigma$ and view the situation from this perspective. Strictly speaking,
 the polarity maps points to hyperplanes, but it defines a duality which maps spaces of dimension $d$ to those of dimension $3-d$ for
 $d=0,\ldots,3$; this is why we speak of the ``polar dual''.  An M-sphere $S$ gives
 rise to a point $S^*$, the pole of the hyperplane carrying the sphere $\Sb \subset \Sigma$.
 An M-circle $k$ in $\R^3$ determines a circle $\kb \subset \Sigma$ in some plane, whose polar image
is a straight line which we denote by $k^*$. Recall that the extended stereographic projection maps the point $S^*$ which determines a sphere $S$ onto the center of $S$. Likewise,
it maps the straight line $k^*$ determining a circle $k$ to the rotational axis of $k$.

A circle family of a Darboux cyclide is completely determined
by a quadratic cone $\Gamma_i$ which carries a family of planes. We now look at the polar dual $\Gamma_i^*$ of the cone,
and as above distinguish the two cases.

\noindent
{\bf Paired circle families.}
If $\Gamma_i$ contains two families of planes and thus determines two paired circle families on a cyclide,
its  polar image $\Gamma_i^*$ is a {\em ruled quadric} in some hyperplane $V_i^*$ (the polar hyperplane
of the cone vertex $V_i$).  Each ruling of $\Gamma_i^*$ is polar to a plane in $\Gamma_i$ and thus determines a circle on the cyclide.
The extended stereographic projection maps the ruled quadric $\Gamma_i^*$ to a quadric (or a plane) $G_i \subset \R^3$.
The axes of the circles in the paired families are projected to two families of lines in $G_i$.

The cone vertex $V_i$ represents an
%(not necessarily real)
M-sphere; let us call it $B_i$ (note: $B_i^* = V_i$) and
consider the associated complex of spheres $S$ which are orthogonal to $B_i$. The poles $S^*$ representing these spheres $S$ lie in the hyperplane $V_i^*$. This implies (recall Remark~\ref{M-symmetry} on M-symmetries) that the cyclide is symmetric with respect to $B_i$. Note that {\it each of the circles in the two families is also symmetric with respect to $B_i$, i.e., orthogonal to $B_i$}, since it is an intersection of two spheres of the complex.

\begin{figure}[htb]
\centering{%\vspace{4cm}
\includegraphics[width=.6\textwidth]{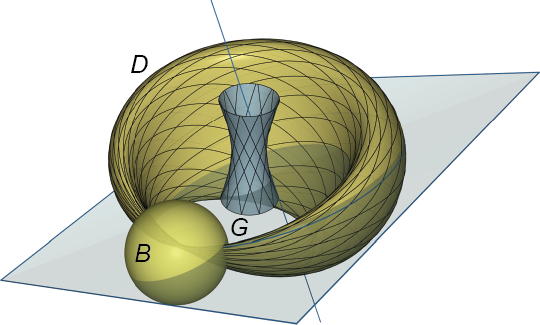}
}
 \caption{A paired family of circles on a cyclide $D$ can be constructed as follows: Prescribe a ruled quadric $G$
 and a sphere $B$, and consider all circles which are orthogonal to $B$ and have a rotational axis on $G$. This
 may be viewed as a transformation which maps the rulings of $G$ to circles on $D$.}
 \label{fig:pairedfamily}
\end{figure}

\begin{prop} \label{prop3} All circles of two paired families are orthogonal to a fixed M-sphere $B_i$ and their axes form the two families of rulings on a ruled quadric $G_i$.  The M-sphere $B_i$ needs not be real and may degenerate into a point; in the latter case, the cyclide is a ruled quadric ($B_i=\infty$) or the image of a ruled quadric under an inversion. The M-sphere $B_i$ may also be plane; in that case, the circle axes are the tangents of a conic in $B_i$ or a pair of pencils of lines in~$B_i$.
\end{prop}

%\mscomm{To be checked by Mikhail: how nonruled quadrics may appear in this construction?}

The latter case appears when the hyperplane
$V_i^*$ passes through the center $Z$ of the stereographic projection and thus is mapped to a plane $B_i$ in $\R^3$.
Now, $G_i$, the image of $\Gamma_i^*$ under the extended stereographic projection, is no longer a ruled quadric,
but degenerates to the union of tangents of a conic or two pencils of lines.

\begin{figure}
\centering{
\includegraphics[width=.7\textwidth]
{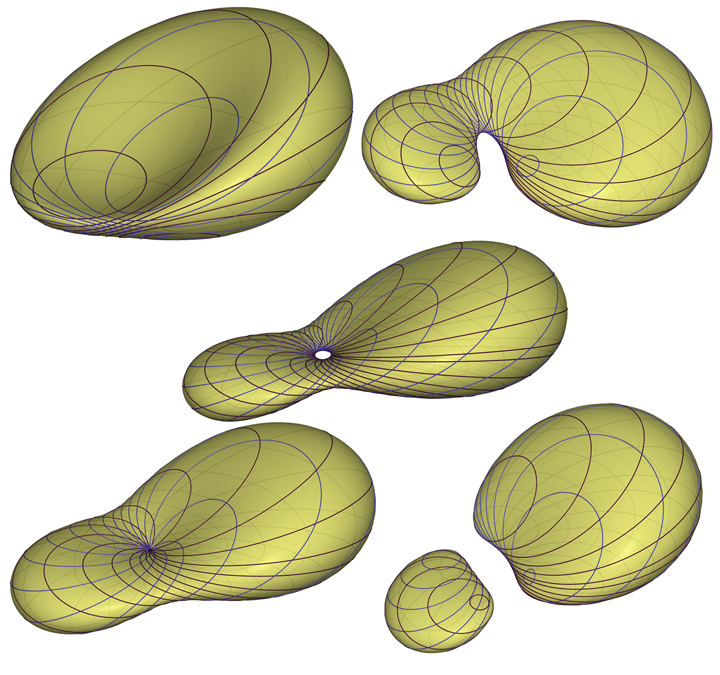}
}
 \caption{Various shapes of Darboux cyclides: A smooth Darboux cyclide is topologically either a sphere (top), or a torus (middle), or a pair of two spheres (bottom right); see \cite{takeuchi:2000}. A cyclide with singularities in the image of a quadric under inversion (bottom left). Paired families of circles on the cyclides are shown.}
\label{fig:manycyclides}
\end{figure}

Thus we can construct paired circle families (and thus also Darboux cyclides) in a very simple way
by prescribing a non-degenerate ruled quadric $G$ and a sphere $B$ (or sphere complex) and proceed according to Proposition~\ref{prop3} (see Fig.~\ref{fig:pairedfamily} and~\ref{fig:manycyclides}).
Note that this construction would fail in the special case where $B$ is a plane: prescribing a conic in a plane
$B$ is not enough to get the cyclide, since circles whose axes are tangent to the conic are automatically
orthogonal to $B$. %{\bf Can we find a direct construction also in this case?}

\begin{rem} For readers familiar with non-Euclidean geometry we mention that Darboux cyclides may
be obtained as follows: Take a ruled quadric in the projective model of hyperbolic or elliptic 3-space
and transfer it (through the well-known Darboux transformation) to the respective conformal model. All
cyclides carrying two paired families of circles which are orthogonal to a real or complex sphere
can be obtained in this way; the former case is that of hyperbolic geometry, the latter belongs to
elliptic geometry.
\end{rem}

\noindent
{\bf Single circle families.}
%-------------------------------------------------------------------------------------------
If the cone $\Gamma_i$ contains only one family of planes and thus has a vertex line $v_i$, its polar image is
a {\it conic} $\Gamma_i^*=:c_i$ in the following sense: The points $S^*$ of the conic $c_i$ determine the spheres $S$ which envelope the
cyclide $D$ and the tangents of $c_i$ represent the circles along which these spheres touch $D$. The plane $P$ of $c_i$ is
the polar dual $v_i^*$ of $v_i$; it defines a bundle of spheres.  Points on the vertex line $v_i$ represent the spheres of a pencil, which
are orthogonal to the spheres of the bundle.
% This shows that the spheres $S$ are orthogonal to the (not necessarily)
%real base circle $b_i$ of the pencil of spheres. In particular, the spheres are centered in the (real) plane of $b_i$.
%The construction of the spheres $S$ enveloping the cyclide $D$ is discussed in section \ref{sec:design}.

Recall that the spheres of a bundle which are orthogonal to an imaginary circle pass through two fixed points; this is the case of Darboux canal surfaces which are quadratic cones (if one of the two points is $\infty$) or images of quadratic cones under an inversion.

\begin{figure}[htb]
\centering{
\includegraphics[width=.5\textwidth]{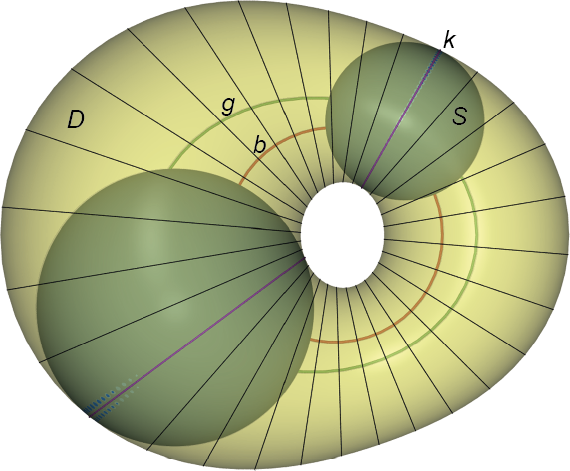}}
 \caption{A cyclide $D$ which carries a single family of circles can be constructed as follows: Prescribe a conic $g$
 and a circle $b$ in the same plane, and consider all spheres $S$ which are centered on $g$ and orthogonal to $b$. The cyclide is the envelope of these spheres. Each sphere $S$ touches the cyclide in a circle $k$ (of the single family), whose axis is a tangent of $g$ and whose plane passes through the axis of $b$.}
 \label{fig:singlefamily}
\end{figure}

\begin{prop} \label{prop4}  A canal surface which is an
irreducible Darboux cyclide (not of revolution)
is the envelope of spheres $S$ which are centered at a conic $g$ and are
orthogonal to an M-circle $b$ in the plane of $g$ (Fig.~\ref{fig:singlefamily}). The M-circle $b$ may degenerate to a point and may also be imaginary; in the latter case, the cyclide is a quadratic cylinder ($b=\infty$), a quadratic cone or an image of a quadratic cone or cylinder under an inversion.
\end{prop}

\begin{rem}
An exceptional case happens if the spheres are all centered on a straight line $g$ (i.e., if the plane of the
conic $c_i$ passes through $Z$), i.e.,  the resulting surface is a rotational surface with
axis $g$. The surface may be defined by rotating a circle around the axis, but only if the surface has another real circle family. There are cyclides for which this is not the case: as an example, take a rotational ellipsoid or
the image of it under an inversion with respect to a co-axial sphere. We are not pursuing the study of cyclides with only one family of real circles, since we are interested in surfaces which carry at least two families of circles.
\end{rem}

%%%%%%%%%%%%%%%%%%%%%%%%%%%%%%%%%%%%%%%%%%%%%%%%%%%%%%%%%
\section{Design tools for Darboux cyclides} \label{sec:design}
%%%%%%%%%%%%%%%%%%%%%%%%%%%%%%%%%%%%%%%%%%%%%%%%%%%%%%%%%

Section \ref{ssec:polardual} shows that we can translate design methods for ruled quadrics and conics to design tools for paired and single circle families on Darboux cyclides. This will be elaborated in Section~\ref{ssec:design} in more detail.
In Section~\ref{ssec:constant-angle} we give an example of circle families intersecting under constant angle and in Section~\ref{ssec:parametrization} we very briefly discuss the parametric representation of Darboux cyclides.

\subsection{Designing circle families on cyclides.} \label{ssec:design}
%---------------------------------------------------------
%
Most of the following material is based on Proposition~\ref{prop3}, which generates paired families with
help of a sphere complex and a ruled quadric. We will now illustrate this at hand of the
case where the sphere complex consists of all spheres orthogonal to a real sphere $B$. It is pretty
straightforward to discuss the other cases.
%(cf. Fig.~\ref{fig:sphere-complexes}).

The following design (input) methods for ruled quadrics will be transferred to cyclide design based on Proposition~\ref{prop3}:
\begin{enumerate}
\item A ruled quadric is defined by three pairwise skew straight lines.
\item A ruled quadric is determined by a skew quad (formed by rulings) and a further point.
\end{enumerate}

\noindent
{\bf (1)} Translating the generation of a ruled quadric $G$ with 3 pairwise skew rulings to
cyclide design is done as follows: We prescribe a sphere $B$ and three circles $k_1, k_2, k_3$ which
are orthogonal to $B$. Since the rulings of $G$ shall be the axes of the circles, we have
to prescribe the circles $k_i$ so that their axes $A_i$ are pairwise skew. Two circles which
are orthogonal to $B$ have skew axes if and only if they do not lie on a common sphere. Having
fixed the circles $k_i$, the cyclide follows easily. One just has to compute a ruled quadric
from 3 pairwise skew lines $A_1, A_2, A_3$. One way is to first compute the rulings of the other
family, as the set of lines intersecting $A_1, A_2, A_3$. To get such a ruling, pick a point $X \in A_1$
and intersect the connecting planes $X \vee A_2$ and $X \vee A_3$.

\begin{figure}[htb]
\centering{
\includegraphics[width=.7\textwidth]{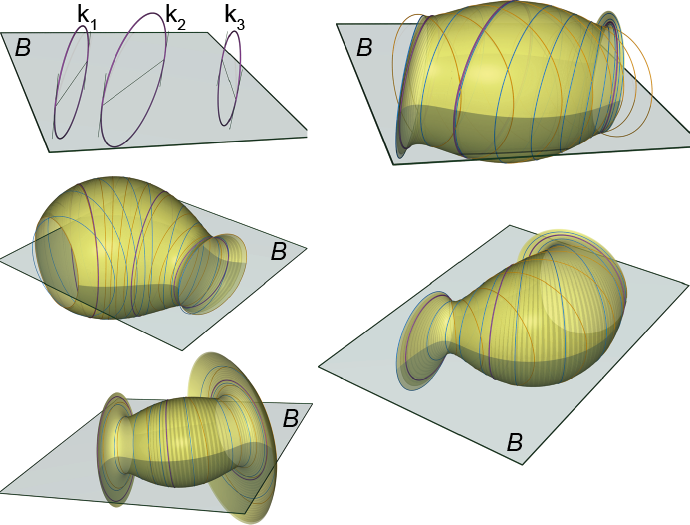}
}
 \caption{Paired families of circles on a Darboux cyclide can be defined by $3$ circles $k_1,k_2,k_3$ which are
 orthogonal to a plane $B$ (and no two of which lie on a common sphere). The figure shows a variety
 of results obtained with this method.}
 \label{fig:cyclide-from3circles}
\end{figure}

It may be more intuitive to let $B$ be a plane and let the user prescribe (half) circles orthogonal
to $B$. Any two such circles should not be co-spherical. Lacking a direct construction in the
case where $B$ is a plane, one can first apply an inversion to map $B$ to a sphere, proceed as above, and
finally transform back with the (same) inversion. We will obtain a cyclide passing through the
three prescribed circles (see Fig.~\ref {fig:cyclide-from3circles}).

\begin{figure}[htb]
\centering{%\vspace{4cm}}
\includegraphics%[width=12.04cm]
{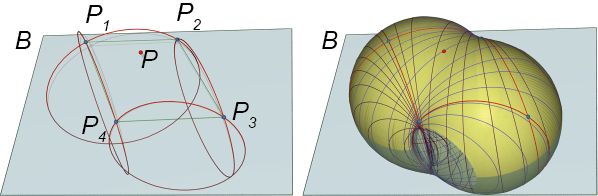}
}
 \caption{A cyclide patch can be defined by its four boundary circles (all being orthogonal to a common
 sphere or plane) and a further point. Especially if $B$ is a horizontal plane, this can be very useful
 for architectural design, since one gets a structure formed by two families of circular arcs in vertical
 planes.}
 \label{fig:cyclide-fromquad}
\end{figure}

{\bf (2)} The second type of input is closely related to the first one, but may be more useful for
practice, since we now prescribe the {\it circular boundaries and a point of
the cyclide patch} (see Fig.~\ref {fig:cyclide-fromquad}).
Let us describe the case with a plane $B$. To this end,
we can input a quad with vertices $P_1,\ldots, P_4$ and a plane $B$. Each edge, say $P_1,P_2$ of the quad can
be reflected at $B$ to obtain points $P_1',P_2'$, and now one has a unique circle passing through
$P_1,P_2, P_1', P_2'$; this circle $k_1$ is orthogonal to $B$. So we end up with four circular arcs, being
admissible boundaries of a cyclide patch. We may then prescribe one more point $P$ of the desired patch
and the cyclide $D$ is determined. This input can be transformed into the above one by constructing a further circle $k$ through $P$: take two opposite circular boundaries, say $k_1, k_3$ and intersect the two connecting spheres $P \vee k_1$ and $P \vee k_3$.

%{\bf Maybe we could make a design based on a collection of cyclide patches.}

%\subsection{Designing Darboux cyclides which are canal surfaces}
%---------------------------------------------------------------
%

\begin{rem}
Proposition \ref{prop4} provides us with a transformation of conics into single families of circles
on cyclides. Thus, in a quite analogous way to the design of paired circle families from ruled quadrics, we can transfer
constructions of conics into those for single circle families on cyclides. This is rather
straightforward, and thus we do not pursue it in the present paper.
\end{rem}

%\textbf{Important: We need to show how to compute the further circle families and give an example.  After all, our focus is on webs.}

\subsection{Two circle families intersecting under a constant angle}\label{ssec:constant-angle}

\begin{figure}[htb]
\centering{
\includegraphics[width=.9\textwidth]{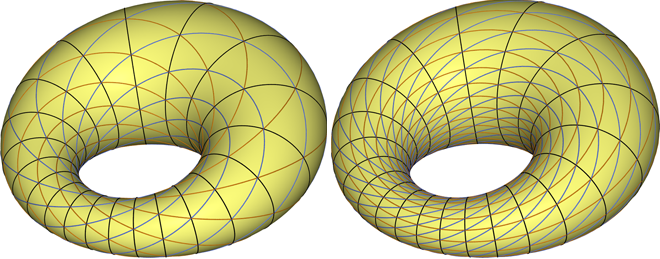}}
\caption{Since the Villarceau circles on a Dupin ring cyclide intersect the principal curvature
 circles under a constant angle, there are several ways of extracting circle families with
 a constant intersection angle. This leads to congruent nodes in an architectural construction.
 We show here two different types of $3$-webs with congruent nodes, as discussed in Theorem \ref{thm-main}.}
 \label{fig:congruent-nodes}
\end{figure}

\cite{ivey-1995} proved that a surface which carries two orthogonal families of circles is a
Dupin cyclide. Apart from the orthogonal network of principal curvature lines (circles) on a Dupin cyclide, there is also the case of the paired families of Villarceau circles on a ring cyclide (lying in doubly tangent planes). Villarceu circles from different families always intersect under a constant angle, as they intersect the principal curvature circles under a constant angle. For specific cyclides this angle can be 90 degrees (see Fig.~\ref{fig:congruent-nodes}).

We see that any pair of circle families on a ring Dupin cyclide exhibits a constant intersection angle, which implies congruent nodes in an architectural construction. We conjecture that the Dupin cyclides are the only surfaces which have two circle families that intersect under a constant angle.

\subsection{Parametric representation of Darboux cyclides} \label{ssec:parametrization}
%----------------------------------
From the point of view of complex algebraic geometry a Darboux cyclide is a Del Pezzo surface of degree $4$ (as a projection of the intersection of two quadrics in $4$-space). In particular, Darboux cyclides are rational surfaces. However, in general they may not have a \emph{real} rational parametrization. For instance, the cyclide $(x^2+y^2+z^2-1)(x^2+y^2+z^2-2)+1/10=0$ is irreducible and disconnected, hence has no proper real rational parametrization.
%Below we introduce a large class of cyclides having a real biquadratic rational parametrization.

For an arbitrary Darboux cyclide $D$ it is not hard to find a real parametrization of the form
$$
x(s,t)=\frac{X_1+X_2\sqrt{P}}{W_1+W_2\sqrt{P}}, \
y(s,t)=\frac{Y_1+Y_2\sqrt{P}}{W_1+W_2\sqrt{P}}, \
z(s,t)=\frac{Z_1+Z_2\sqrt{P}}{W_1+W_2\sqrt{P}},
$$
where $X_1,X_2,Y_1,Y_2,Z_1,Z_2,W_1,W_2,P$ are polynomials in $s$ and $t$ of degree at most $4$.

Indeed, let ${\mathbf X}^{{T}}\cdot A \cdot {\mathbf X}=0$ be the equation of~$\bar D$. Since the polynomial $\det (A-tJ)$ has degree $5$, it has a real root $t_i$. Thus $\bar D$ is the intersection of the cone ${\mathbf X}^{{T}}\cdot (A-t_iJ)\cdot {\mathbf X}=0$ and the sphere $\Sigma$. Depending on the position of the vertex $V_i$ of the cone with respect to the sphere $\Sigma$ there are the following $3$ cases.

{\it Case 1}: $V_i$ is in the ``surface'' of the sphere $\Sigma$. Consider the hyperplane $H_i$ touching the sphere $\Sigma$ at the point opposite to $V_i$. The hyperplane $H_i$ intersects the cone by a quadric $D_i$. Take a quadratic parametrization of the quadric $D_i$. Composing it with the inverse stereographic projection from the vertex $D_i$ and then with the stereographic projection from the point $Z$ we get a rational parametrization of the cyclide $D$ of degree a most $4$.

{\it Case 2}: $V_i$ is outside the sphere $\Sigma$. Perform a projective transformation $\kappa$ of the $4$-space preserving $\Sigma$ and taking $V_i$ to the infinitely distant point in the direction of the $x_1$-axis. The image of the cone under the transformation $\kappa$ is a cylinder crossing
the hyperplane $x_1=0$ by a quadric (or a plane) $D_i$. Thus the quadric (or a plane) $D_i$ is the projection of $\kappa \bar D$ along the $x_1$-axis. Take a quadratic parametrization of the quadric $D_i$. Inverting the projection, performing the inverse transformation $\kappa^{-1}$, and finally performing the stereographic projection from $Z$ we get the required parametrization.

{\it Case 3}: $V_i$ is inside the sphere $\Sigma$. Perform a projective transformation $\kappa$ of the $4$-space preserving $\Sigma$ and taking $V$ to the center of the sphere $\Sigma$. The image of the cone under the transformation intersects the hyperplane $x_1=1$ by a quadric (or a plane) $D_i$. Thus the quadric (or a plane) $D_i$ is the central projection of the cyclide from the origin. Take a quadratic parametrization of the quadric $D_i$. Inverting the projection, performing the inverse transformation $\kappa^{-1}$, and finally performing the stereographic projection from $Z$ we get the required parametrization.

%%%%%%%%%%%%%%%%%%%%%%%%%%%%%%%%%%%%%%%%%%%%%%%%%%%%%%%%%%
\section{Webs from circles on Darboux cyclides} \label{sec:cwebs}
%%%%%%%%%%%%%%%%%%%%%%%%%%%%%%%%%%%%%%%%%%%%%%%%%%%%%%%%%

%{\bf Definition of a web required: discrete and continuous webs are often confused in the paper. Probably, a closing condition should be stated, methods of projections and symmetry briefly explained, references given. Convention: by default all webs are hexagonal. --- To be done by Mikhail.}

In this section we define $3$-webs, give classical examples of planar C-webs, and finally classify all possible C-webs on Darboux cyclides.

\begin{figure}
\centering{
\includegraphics[width=\textwidth]{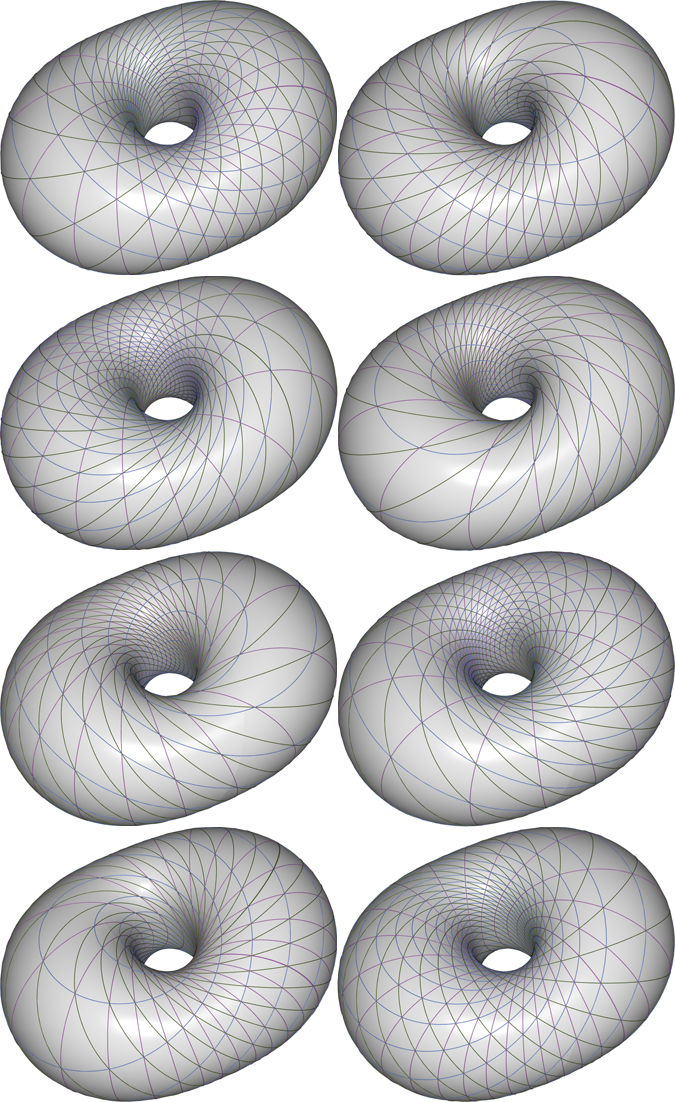}}
\caption{Eight different $3$-webs of type 1 in a Darboux cyclide.}
\label{fig:closed-web}
\end{figure}

A \emph{$3$-web} (or \emph{hexagonal web}) in a surface is $3$ families of smooth curves which are locally diffeomorphic to $3$ families of lines $x=\mathrm{const}$, $y=\mathrm{const}$, $x+y=\mathrm{const}$ in the plane $Oxy$. These lines for integral values of the constants form a triangular lattice in the plane; taking the preimage of the lattice under the diffeomorphism, we get a \emph{discrete $3$-web} in the surface (see Fig.~\ref{fig:closed-web}).

\begin{figure}
\centering{
\includegraphics[width=.6\textwidth]{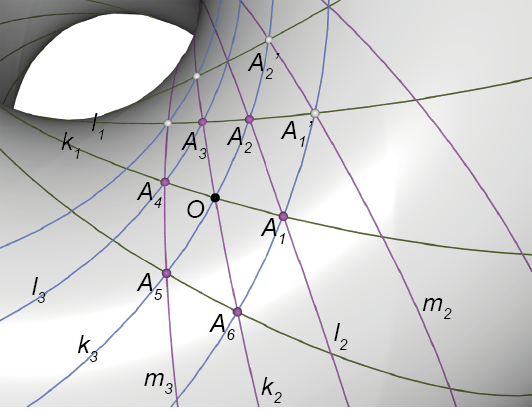}
}
\caption{The closure condition characterizing $3$-webs; see the beginning of Section~\ref{sec:cwebs} for explanation. This figure also illustrates an algorithm for computation of discrete $3$-webs on Darboux cyclides; see the end of Section~\ref{sec:cwebs}.}
\label{fig:finite-web}
\end{figure}

A {$3$-web} is characterized by the following \emph{closure condition} (\cite{blaschke-bol:1938}; see also Fig.~\ref{fig:finite-web}). Take an arbitrary point $O$ of the surface. Draw the curves $k_1$, $k_2$, $k_3$ of the $3$ families through the point $O$. Take an arbitrary point $A_1\in k_1$ sufficiently close to~$O$. Draw the curve $l_2$ of the second family through the point $A_1$. Set $A_2=l_2\cap k_3$. Draw the curve $l_1$ of the remaining (first) family through the point $A_2$. Set $A_3=l_1\cap k_2$. Draw the curve $l_3$ of the remaining (third) family through~$A_3$ and set $A_4=l_3\cap k_1$. Continuing in this way, construct points $A_5$, $A_6$, and $A_7$. The \emph{closure condition} asserts that $A_7=A_1$. This condition and symmetry considerations often allow to verify that certain triple of families is indeed a $3$-web.

First let us give the following known beautiful examples of $3$-webs of lines and circles in the plane needed for the sequel. In these examples by a \emph{$3$-web} we mean a  $3$-web in an appropriate nonempty domain in the plane. Three families of curves are \emph{crossing}, if some $3$ curves of distinct families have a common point of transversal intersection.

\begin{thm}\label{thm-Graf-Sauer} \textup{(\cite{graf:sauer:1924, blaschke-bol:1938})}
Three crossing families of (pairwise distinct) straight lines in the plane form a $3$-web if and only if all these lines are tangent to one curve of algebraic class $3$ (including reducible ones).
\end{thm}

A %\emph{bicircular quartic}
\emph{cyclic}
is a curve given by equation of the form
$$
\lambda(x^2+y^2)^2+(x^2+y^2)(\mu x+\nu y)+Q(x,y)=0,
$$
where $\lambda,\mu,\nu\in\mathbb{R}$ and $Q(x,y)$ is a polynomial of degree at most $2$. (This is a one-dimensional counterpart of a cyclide.)

\begin{thm} \label{thm-Wunderlich} \textup{(\cite{wunderlich-38})} Three crossing families of circles doubly tangent (in the sense of complex %projective
algebraic geometry) to a %bicircular quadric
cyclic (possibly with singularities)
form a $3$-web.
\end{thm}

\begin{rem} \label{rmk-Wunderlich} Performing the inverse stereographic projection and then a projective transformation one gets the following corollary:
%result (from the Wunderlich theorem):
the tangent planes to $3$ of the $4$ quadratic cones (possibly reducible) passing through the intersection $c$ of two oval quadrics in $\mathbb{R}^3$ cut a $3$-web of conics in each of the oval quadrics, once the $3$ families of conics are crossing (cf.~Fig.~\ref{fig:quadpencil}).
\end{rem}

Since we have up to $6$ real families of circles on a cyclide, the question arises whether and how
we can arrange triples of them into  $3$-webs. To this
end, we will prove the following result (see Fig.~\ref{fig:closed-web}):

\begin{thm} \label{thm-main}  Three families of circles on a nontrivial irreducible Darboux cyclide
%(different from a sphere or plane)
form a  $3$-web unless one takes two nonspecial paired families and another family which has a paired one. Thus we have $5$ types of $3$-webs from circles on a Darboux cyclide:
\begin{enumerate}
    \item[(i)] A web of type 1 is formed by $3$ non-single families of circles, such that no two of them are paired families. Hence, this type only exists on a cyclide with $6$ real families of circles, and we have $8$ different $3$-webs on such a cyclide; see Fig.~\ref{fig:closed-web}.
    \item[(ii)] A web of type 2 is formed by $2$
    special paired families and another family which has a paired one. This type only exists on the image of a one-sheet hyperboloid under an inversion; see Fig.~\ref{fig:web-type2}.
    \item[(iii)] A web of type 3 is formed by a single family and $2$ paired families; see Fig.~\ref{fig:congruent-nodes} to the left.
    \item[(iv)] A web of type 4 is formed by a single family and $2$ non-single families, which are not two paired families; see Fig.~\ref{fig:web-type4}.
    \item[(v)] A web of type 5 consists of $2$ single families and another family which has a paired one. Hence, this type only exists on the Dupin ring cyclide and consists of its principal curvature lines and a family of Villarceau circles; see Fig.~\ref{fig:congruent-nodes} to the right.
\end{enumerate}
\end{thm}

\begin{proof} The main idea of the proof is to use the  surface $\Db \subset \R^4$, the cones which
define the circle families and to verify webs
at hand of appropriate projections onto a plane or 3-space. Types 2,4,5 are in fact degenerations of type 1.
%The most nontrivial is type 1. All the other types can be obtained by various degenerations.

\noindent {\it Type 1:} Here $\Db$ is the carrier of a pencil of quadrics which contains three
cones (with vertices $V_1, V_2, V_3$) each of which contains %two families of
planes. By Remark~\ref{rem:largest eigenvalue} the pencil of quadrics must also contain a cone with a real vertex $V_4$ and an oval base quadric. Now we project $\Db$ from $V_4$ onto a hyperplane $H_4$ (not through
$V_4$) and obtain (a part of) an oval quadric $\Phi$. The three circle families  of $\Db$ (no two of them paired) are projected to conics (on $\Phi$) in tangent planes of quadratic cones; $\Phi$ and these cones lie in a pencil of quadrics. By Proposition~\ref{prop-1commonpoint} it follows that the conics of distinct families are crossing.
We have now %(up to a projective transformation)
exactly the situation of Remark~\ref{rmk-Wunderlich} and conclude that these three families of conics form a $3$-web. Hence, also the $3$ circle families on the cyclide form a $3$-web.

\noindent {\it Type 2:} Let $V_1\in\Sigma$ be vertex of the cone giving the special paired families. Let $V_2$ be the vertex of the cone $\Gamma_2$ giving the other non-single family. Project from the line connecting $V_1$ and $V_2$ onto a plane. Then the special paired families are projected either to tangents to a conic or to two pencils of lines. Since $V_1\in\Sigma$ it follows that $V_1\in\Gamma_2$. Thus the 3rd non-single family is projected to a pencil of lines. The resulting straight lines are tangent to a reducible curve of algebraic class $3$ and by Theorem~\ref{thm-Graf-Sauer} form a $3$-web.
Hence, also the $3$ circle families on the cyclide form a $3$-web.

Special paired families only exist on cyclides which are images of quadrics under inversion. Among the quadrics only a one-sheet hyperboloid (not of revolution) contains $4$ families of circles/lines. Thus the cyclide in type 2 can only be an image of a one-sheet hyperboloid under inversion. (This web on a hyperboloid was also discussed by \cite{wunderlich-38}.)

\begin{figure}
\centering{
\includegraphics[width=0.5\textwidth]{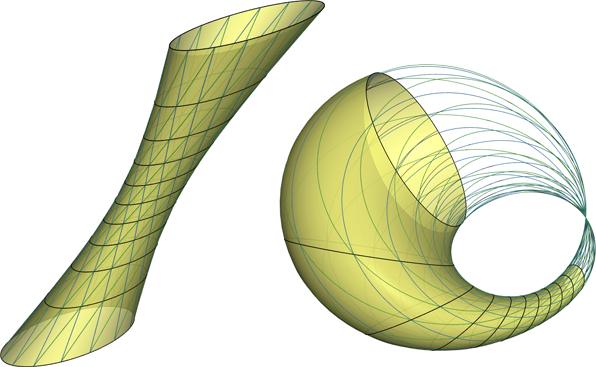}}
\caption{A $3$-web of type 2 in a Darboux cyclide (right) obtained from a non-rotational one-sheet  hyperboloid  (left) by an inversion.}
\label{fig:web-type2}
\end{figure}

\noindent {\it Type 3:} %\mscomm{Needs more attention.}
We connect the cone vertex $V_1$ (associated with the two paired families) with any
point on the vertex line $v_2$ of the cone that determines the single family. From the resulting line $L$
we project into a plane (not intersecting $L$). The image of the two paired circle families are either the tangents
of a conic or two pencils of lines. The image of the single family is a pencil of lines. So together we have the tangents of a reducible curve of algebraic class $3$, which according to Theorem~\ref{thm-Graf-Sauer} form a $3$-web.

\noindent {\it Type 4:} The proof is literally the same as in the proof for type 1. We only need to remark that one of the cones in the hyperplane $H_4$ now degenerates to a pair of planes and its \emph{tangent planes} are understood in the sense of algebraic geometry: these are the ones passing through the vertex line of the cone.

\begin{figure}
\centering{
\includegraphics[width=.6\textwidth]{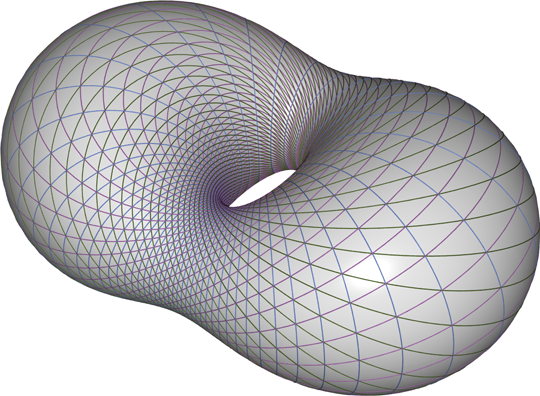}
}
\caption{A $3$-web of type 4 in a Darboux cyclide.}
\label{fig:web-type4}
\end{figure}

\noindent {\it Type 5:}
We may apply an inversion to map the ring Dupin cyclide to a ring torus. Now it follows from the rotational symmetry that the Villarceau circles together with meridian circles
and parallel circles (rotational paths) form a $3$-web.

We still have to show that two (nonspecial) paired families and a 3rd (non-single) circle family cannot form a $3$-web. Again, project from the line connecting the two involved cone vertices, say $V_1, V_2$ onto a plane. The two paired families are projected either to tangents to a conic or to two pencils of lines. Since $V_1\not\in\Sigma$ it follows that the 3rd family is projected to tangents to a (nondegenerate) conic. The resulting straight lines are tangent to a reducible curve of algebraic class $4$ and
due to Theorem~\ref{thm-Graf-Sauer} do not form a   $3$-web.
\end{proof}

\begin{rem} %One can see easily that the theorem \mscomm{(without the characterization of surfaces having each type of webs)} and its proof remain true if on replaces ``\mscomm{nontrivial irreducible} Darboux cyclide'' by ``intersection of two \mscomm{$3$-dimensional} nonsingular quadrics in $\mathbb{R}^4$ \mscomm{containing no $2$-dimensional quadrics as a whole}'' and ``circles'' by ``conics''.
Similarly, there are the same $5$ types of $3$-webs of conics in the intersection of two arbitrary $3$-dimensional nonsingular quadrics in $\mathbb{R}^4$ (if the intersection does not  contain  $2$-dimensional quadrics as a whole).
If we project such a $3$-web from
$\mathbb{R}^4$ to $\mathbb{R}^2$ we get a $3$-web of conics in the plane. The conics of the planar $3$-web are tangent to the ``visual contour'' (i.e., the boundary of the projection) of the intersection of the quadrics, which is a curve of order $\le 12$; cf.~\cite{schmid:1936}, \cite{timorin:2007}.
\end{rem}

%\textbf{Unsure if the latter result is new because there was a lot of papers on webs from conics. Definitely this construction was known for the particular case of projection of Dupin cyclides from $R^3$ to $R^2$.}

%{\bf Still to be done:
%\begin{itemize}
%\item briefly describe the actual computation of webs (with a simple explanatory figure) and the optimization leading to closed webs.
%\end{itemize}
%}

The closure condition can be used for  the computation of
discrete $3$-webs as follows.
%The $3$-webs of the whole cyclide can be generated by the hexagonal web.
Consider $3$ families of curves on a surface, forming a $3$-web. We use the notation from the definition of the closure condition at the beginning of Section~\ref{sec:cwebs}.
We prescribe the requested ``size'' $N$ of the web, as well as starting points $O$ and $A_1\in k_1$. One step of algorithm is to draw the curve $m_3$ of the third family through the point $A_1$; set $A_1'=m_3 \cap l_1$; then draw the curve $m_2$ of the second family through the point $A_1'$; set $A_2'=m_2 \cap k_3$. The next step of algorithm is the same, but now $A_2$ and $A_1'$ are chosen as starting points. Performing $N$ such steps one gets $N$ curves of the first family and $N$ curves of the second one. After that, we make similar $N$ steps starting from the points $O$ and $A_2$ to get $N$ curves of the third family.

One can apply this algorithm to construct \emph{finite} $3$-webs, i.e., discrete $3$-webs forming a triangulation of the whole closed surface.
%(this is only possible when the surface has ring shape).
Indeed, parameterize one of the curve families, say, the first one, so that the curve $k_1$ corresponds to  parameter values both $0$ and $1$ (we assume that the surface has ring shape). Running the algorithm we get $N$ curves of the first family, let $\sigma$ be the parameter value of the $N$-th curve. Let $\tilde A_1\in k_1$ be a point such that $O\tilde A_1=OA_1/\sigma$.
Repeat the same procedure starting from the points $O$ and $\tilde A_1$ until $\sigma$ is close enough to $1$. Eventually a discrete $3$-web close to a finite one is constructed this way.
On Darboux cyclides, there are arbitrary dense finite $3$-webs of type 1 by \cite[Section~4]{wunderlich-38}.
%Figures in this section provide evidence of finite C-webs on Darboux cyclides although we have no proof of their existence so far.

%%%%%%%%%%%%%%%%%%%%%%%%%%%%%%%%%%%%%%%%%%%%%%%%%%%%%%%%%
 \section{Related studies and open problems} \label{sec:future}
%%%%%%%%%%%%%%%%%%%%%%%%%%%%%%%%%%%%%%%%%%%%%%%%%%%%%%%%%

There are several directions in which we would like to continue our research:

 \begin{itemize}
 \item We are still lacking a proof for the conjecture that any surface which carries three families of circles is a Darboux cyclide and that any surface with two families of circles intersecting under a constant angle is a Dupin cyclide. Moreover, we would like to classify all surfaces which carry two families of circles and
 discuss them from a geometric modeling perspective.
 \item Cyclides of isotropic geometry may also deserve attention. Recall that an \emph{isotropic circle} is either a parabola with a vertical axis (a favorite element of architecture) or an ellipse whose projection into a horizontal plane (top view) is a circle. An \emph{isotropic cyclide} is given by equation~\eqref{D}, in which both instances of $(x^2+y^2+z^2)$ are replaced by $(x^2+y^2)$.
     Viewing isotropic cyclides from the Laguerre geometric perspective, we can transform them into surfaces which are enveloped by several families of right circular cones. This would, for example, provide access to architectural structures from strips of right circular cones.
 \item It is also interesting to obtain an algorithm for rational parametrizations of Darboux cyclides and the corresponding rational Bezier representation; this is expected to result from the work of \cite{krasauskas:2011}.
 \end{itemize}

\subsection*{Acknowlegdements}
The authors are grateful to A.~Kochubei for valuable discussions.

\bibliographystyle{elsart-harv}
\bibliography{webs}
\end{document}